\documentstyle[12pt]{article}
\def\C{C\!\!\!\!I}
\def\D{{\cal D}}

\def\E{{\cal E}}
\def\EE{{\bf E}}
\def\F{{\cal F}}
\def\FF{{\bf F}}
\def\G{{\cal G}}
\def\GG{{\bf G}}

\def\K{{\cal K}}
\def\lra{\longrightarrow}
\def\L{{\cal L}}
\def\M{{\cal M}}

\def\O{{\cal O}}
\def\P{{\cal P}}

\def\rh{{\cal R\!\!H}}

\def\t{\tau}

\def\un{\underline}
\def\U{{\cal U}}
\def\V{{\cal V}}
\def\W{{\cal W}}

\def\Z{Z\!\!\!Z}

\def\spec{\mathop{\rm Spec}\nolimits}
\let\ov\overline

\newtheorem{theorem}{Theorem}[section]

\newtheorem{definition}[theorem]{Definition}
\newtheorem{example}[theorem]{Example}
\newtheorem{lemma}[theorem]{Lemma}
\newtheorem{proposition}[theorem]{Proposition}
\def\rem{\refstepcounter{theorem}\paragraph{Remark \thetheorem}}
\def\rems{\refstepcounter{theorem}\paragraph{Remarks \thetheorem}}
\def\proof{\paragraph{Proof}}

\makeatletter \def\l@section{\@dottedtocline{1}{0em}{1.2em}} \makeatother

\begin{document}

\title{Moduli of regular holonomic $\D$-modules \\
with normal crossing singularities }
\author{Nitin Nitsure}
 
\date{09-III-1997}
\maketitle

\centerline{Tata Institute of Fundamental Research, Mumbai 400
005, India.} 
\centerline{e-mail: nitsure@math.tifr.res.in}

\centerline{Mathematics Subject Classification: 14D20, 14F10,
32G81, 32C38} 

\begin{abstract} 
This paper solves the global moduli problem for regular
holonomic $\cal D$-modules with normal crossing singularities 
on a nonsingular complex projective variety. This is done by
introducing a level structure (which gives rise to ``pre-$\cal
D$-modules''), and then introducing a notion of (semi-)stability
and applying Geometric Invariant Theory to construct a coarse
moduli scheme for semistable pre-$\cal D$-modules. A moduli is
constructed also for the corresponding perverse sheaves, and the
Riemann-Hilbert correspondence is represented by an analytic
morphism between these moduli spaces.
\end{abstract}

\vfill
\newpage

\tableofcontents














\vfill
\newpage

\section{Introduction} 
The moduli problem for regular holonomic $\D$ modules on a
non-singular complex projective variety $X$ has the following 
history. Around 1989, Carlos Simpson solved the problem in the
case when the $\D$-modules are $\O$-coherent, which first
appeared in a preliminary version of his famous paper [S]. 
In this case, a $\D$-module $M$ on $X$ is the same as a vector
bundle together with an integrable connection. 
The next case, that of meromorphic connections with regular
singularity along a fixed normal crossing divisor $Y\subset X$
was solved in [N]. As explained there, one has to consider a
level structure in the form of logarithmic lattices for the
meromorphic connections in order to
have a good moduli problem (or an Artin algebraic stack), and
secondly, a notion of semi-stability has to be introduced in
order to be able to apply Geometric Invariant Theory. 
In collaboration with Claude Sabbah, a more general case was
treated in [N-S], where the divisor $Y\subset X$ is required
to be smooth, but the only restriction on the regular holonomic
module $M$ is that its characteristic variety should be
contained in $X\cup N^*_{S,X}$ (this is more general than being a
non-singular or meromorphic connection). Here, we introduced the
notion of pre-$\D$-modules, which play the same role for these
regular holonomic modules that logarithmic connections play for 
regular meromorphic connections. A notion of semistability was
introduced for the pre-$\D$-modules, and a moduli was constructed.
We also constructed a moduli for the corresponding perverse
sheaves, and showed that the Riemann-Hilbert correspondence
defines an analytic morphism from the first moduli to the
second, and has various good properties. 
This is already the most general case if $X$
is $1$ dimensional.

The present paper solves the moduli problem in the general case 
where we have a divisor with {\sl normal crossings}, and the
caracteristic variety of the regular holonomic $\D$-modules
is allowed to be any subset of the union of the conormal bundles
of the nonsingular strata of the divisor. 
This is done by extending the notion of pre-$\D$-modules to this
more general case, defining semistability, and constructing a
moduli for these using GIT methods and Simpson's construction in
[S] of moduli for semistable $\Lambda$-modules.  
Also, a moduli is constructed for the corresponding
perverse sheaves, and the Riemann-Hilbert correspondence is
represented by an analytic morphism having various good
properties. We now give a quick overview of the contents of this
paper.

Let $X$ be a nonsingular variety of dimension $d$, and let
$Y\subset X$ be a divisor with normal crossings (the irreducible
components of $Y$ can be singular). Let $S_d=X$, 
$S_{d-1}=Y$, 
and let $S_i$ be the singular locus of $S_{i+1}$
for $i< d-1$.  This defines a filtration of $X$ by closed
reduced subschemes 
$$S_d\supset S_{d-1}\supset \ldots \supset S_0$$ 
where each $S_i$ is either empty or of pure dimension $i$.
Let $S'_d=X-Y$ and for $i<d$ let $S'_i =S_i-S_{i-1}$. 
The stratification $X= \cup _iS'_i$ is called the {\bf singularity
stratification} of $X$ induced by $Y$. 
Let $T^*X$ be the total space of the cotangent bundle of $X$,
and for $i\le d$ let $N^*_i\subset T^*X$ be 
the locally closed subset which is the conormal bundle
of the closed submanifold $S'_i$ of $X-S_{i-1}$. In particular,
$N^*_d$ is the zero section $X\subset T^*X$. 
Let $N^*(Y)\subset T^*X$ be defined to be the union 
$$N^*(Y) = N^*_d \cup N^*_{d-1} \cup N^*_{d-2} \cup \ldots$$
of all the $N^*_i$ for $i\le d$. 
Note that $N^*(Y)$ is a closed lagrangian subset of $T^*X$, and
any irreducible component of $N^*(Y)$ is contained in the closure
$\ov{N^*_i}$ for some $i$. 

In this paper we consider those regular holonomic $\D$ modules
on $X$ whose characteristic variety is contained in $N^*(Y)$.
Equivalently (under the Riemann-Hilbert correspondence), we
consider perverse sheaves on $X$ which are cohomologically constructible 
with respect to the singularity stratification 
$$X=\cup_{0\le i\le d}S'_i$$ 
Such regular holonomic $\D$ modules (equivalently, such
perverse sheaves) form an abelian category, in which each object
is of finite length. These will be called
regular holonomic $\D$-modules (or perverse sheaves) on $(X,Y)$.

In section 2, we introduce the basic notation involving various
morphisms which arise out of the singularity strata, their
normalizations, and their \'etale coverings.

One of the problems we had to overcome was to give a convenient
$\O$-coherent description of regular holonomic $\D$-modules on
$(X,Y)$. This is done in section 3, by extending the notion of a
pre-$\D$-modules from [N-S] to this more general setup.  

In section 4, we explain the functorial passage from
pre-$\D$-modules to $\D$-modules. In the case of a smooth divisor,
this was only indirectly done, via Malgrange's
presentation of $\D$-modules, in [N-S]. Here we do it more
explicitly in our more general situation. 

In section 5, we define a notion of
(semi-)stability, and construct a moduli scheme for
(semi-)stable pre-$\D$-modules on $(X,Y)$ with prescribed 
numerical data (in the form of Hilbert polynomials). 
We have given here an improved quotient construction, which
allows us to give a much more simplified treatment of stability
when compared with [N-S].

In section 6, we represent perverse sheaves on $(X,Y)$ with 
prescribed numerical data
in {\sl finite terms}, via a notion of `Verdier objects', generalizing
the notion of Verdier objects in [N-S]. 
We then construct a moduli scheme for these Verdier objects.
This is a quotient of an affine scheme by a reductive group, so
does not need any GIT. Though a more general moduli construction
due to Gelfand, MacPherson and Vilonen exists in literature (see
[G-M-V]), the construction here is particularly suited for the study
the Riemann-Hilbert morphism (in section 7), which represents
the de Rham functor.   

In section 7 we show that the Riemann-Hilbert correspondence
defines an analytic morphism from (an open set of) the moduli of
pre-$\D$-modules to the moduli of perverse sheaves. We extend
the rigidity results in [N] and [N-S] to this more general situation,
which in particular means that this morphism is a local
isomorphism at points with {\bf good residual eigenvalues} as
defined later.

{\footnotesize
{\bf Erratum} I point out here some mistakes which have remained
in [N] and [N-S], and their corrections.

(1) The lemma 2.9 of [N] is false as stated, and needs the
additional hypothesis that $A$ is reduced. This additional
hypothesis of reducedness is satisfied in the part of
proposition 2.8 of [N] where this lemma is employed.

(2) On page 58 of [N-S], in the construction of a local
universal family, we take the action of $G_i=PGL(p_i(N))$ on
$Q_i$. This should read $G_i=SL(p_i(N))$ and not $PGL(p_i(N))$. 

(3) In Theorem 4.19.3 on page 63 of [N-S], the statement
``The S-equivalence class of a semistable reduced module $\EE $
equals its isomorphism class if and only if $\EE $ is stable''
should be corrected to read ``The S-equivalence class of a
semistable reduced module $\EE$ equals its isomorphism class if
$\EE$ is stable''. The ``... and only if...'' should be removed.
(The proof only proves the ``if'' part, ignores the ``only if''
part, and the ``only if'' statement is in fact trivially false.) 

}

{\bf Acknowledgement} I thank the International Center for
Theoretical Physics, \\ 
Trieste, for its hospitality while part of
this work was done. I thank H\'el\`ene Esnault for providing the
interesting example \ref{Esnault} below, in answer to a question.

\section{Preliminaries on normal crossing divisors}
In section 2.1 we define various objects naturally associated
with a normal crossing divisor. The notation introduced here is
summarized in section 2.2 for easy reference, and is used
without further comment in the rest of the paper.

\subsection{Basic definitions}

Let $X$ be a nonsingular variety of dimension $d$, $Y$ a normal
crossing divisor, and let closed subsets $S_i\subset X$ for
$i\le d$ be defined as in the introduction ($S_d=X$,
$S_{d-1}=Y$, and by descending induction, $S_{i-1}$ is the
singularity locus of $S_i$). Let $m$ be the smallest integer
such that $S_m$ is nonempty. 
We allow the possibility that $m$ is any integer from $0$ to $d$
(for example, if $m=d$ then $Y$ is empty, and if $m=d-1$ 
when $Y$ is smooth).

For $i\ge m$, let $X_i$ be the normalization of $S_i$, with
projection $p_i:X_i\to S_i$.
This is a reduced, nonsingular, $i$-dimensional scheme of finite
type over $\C$ for $i\ge m$. This may not be connected (same as
irreducible), and we denote its components by 
$X_{i,a}$, as $a$ varies over the indexing set $\pi_0(X_i)$.
Note that the fiber of $p_i:X_i\to S_i$ over a point $y\in S_i$ is
the set of all branches of $S_i$ which pass through $y$ (where
by definition a branch means a component in the completion of
the local ring).   
Now let $I$ be a nonempty subset of $\{\, m,\ldots ,d-1\,\}$,
and for any such $I$, let $m(I)$ denote the smallest element of $I$.
To any such $I$, we now associate a scheme $X_I$ and a morphism  
$p_I:X_I\to S_{m(I)}$ as follows. By definition, 
$X_I$ is the finite scheme over $S_{m(I)}$ whose fiber over 
any point $y$ of $S_{m(I)}$ consists of all nested sequences $(x_j)$
of branches of $S_j$ at the point $y$, for $j$ varying over $I$.
`Nested' means for any two $j,\,k\in I$ with $j \le k$, $x_k\in
X_k$ is a 
branch of $S_k$ containing the branch $x_j\in X_j$ of $S_j$. 
In particular when $I=\{\, i\,\}$,   
$p_{\{ i\}}:X_{\{ i\} }\to S_i$ is just the normalization 
$p_i:X_i\to S_i$ of $S_i$. 

For nonempty subsets $I\subset J\subset \{\, m,\ldots ,d-1\,\}$,
we have a canonical forgetful map $p_{I,J}:X_J\to X_I$. If $I\subset
J\subset K$, then by definition we have the equality
$$p_{I,K}=p_{I,J}\circ p_{J,K}:X_K\to X_I$$
For $m\le i\le d-1$ we denote by 
$Y^*_i$ the $d-i$ sheeted finite \'etale cover $p_{\{ i\},\{ i+1
\} }: X_{\{ i,i+1 \} }\to X_i$ of $X_i$, which splits the branches of
$S_{i+1}$ which meet along $S_i$. 
Similarly, we denote by $Z_i$ the $C^{d-i}_2$ sheeted finite
etale cover $p_{\{ i\} ,\{ i+2\} }: X_{\{i,i+2\} } \to X_i$ of $X_i$,
which splits the branches of $S_{i+2}$ which meet along $S_i$.
We denote by $Z^*_i$ the $2$ sheeted finite \'etale cover
$p_{\{ i, i+2\} ,\{ i,i+1,i+2\} }: X_{\{i,i+1,i+2 \} } \to X_{\{
i,i+2 \} }$ of $Z_i$. Let $f_i:X_i\to X$ be the composite
$X_i\stackrel{p_i}{\to}S_i\to X$. 

Let $N_i$ be the normal bundle of $X_i$ in $X$.
This is defined by means of the exact sequence
$$ 0\to T_{X_i} \to f_i^*T_X \to N_i\to 0 $$
Let $T(Y)\subset T_X$ be the closed subset of $T_X$ consisting of
vectors tangent to branches of $Y$. This gives a closed subset
$f_i^*T(Y)\subset f_i^*T_X$. Let $F_i$ be the closed subset of
$N_i$ which is the image of $f^*T(Y)$ under the morphism
$f^*T_X\to N_i$ of geometric vector bundles. 
Similarly, let $N_{i,i+1}$ be the normal bundle to
$Y^*_i=X_{i,i+1}$ in $X$, defined with respect to the composite
morphism 
$Y^*_i\to X_i\to X$, and let $F_{i,i+1}\subset N_{i,i+1}$ 
be the normal crossing divisor in the total space of
$N_{i,i+1}$, defined by vectors tangent to branches of $Y$. 

Let $Y^*=Y^*_{d-1} =X_{d-1}$ be the normalization of
$Y$. Let $h_i:Y^*_i\to Y^*$ be 
the canonical map, which associates to a point $(x_i, x_{i+1})\in
X_{\{i,i+1\} }= Y^*_i$ the unique branch $y_{d-1}\in Y^*$ of $Y$ at 
$p_i(x_i)\in S_i$ such that 
$$x_i = x_{i+1}\cap y_{d-1}$$
This defines a vector hyper subbundle $H_{i,i+1} \subset
N_{i,i+1}$ (hyper subbundle means rank is less by $1$) of
$N_{i,i+1}$, which is given by vectors tangent to the branch of
$Y$ given by $h_i:Y^*_i\to Y^*$. Note that $H_{i,i+1}$ is a
nonsingular irreducible component of $F_{i,i+1}$. 

We now define some sheaves $\D_i$ of algebras 
differential operators on $X_i$ and $\D^*_i$ on $Y^*_i$. 
(These are `split almost polynomial algebras' of differential
operators in the terminology of Simpson [S] as explained in
later in section 5.1.) 
For this we need the following general remark: 

\rem Let $V$ be nonsingular, $M\subset V$ a divisor with normal
crossing, and $M'\subset M$ a nonsingular component of $M$, with
inclusion $f:M'\hookrightarrow V$. Let $\Lambda = \D_V[\log M]$
be the subring of $\D_V$ which preserves the ideal sheaf
$I_M\subset \O_V$, and let
$\Lambda'= \O_{M'}\otimes_{\O_V}\Lambda$ (we will write
$\Lambda' =f^*(\Lambda)$ or $\Lambda'=\Lambda|M'$ for brevity,
though this abbreviated notation conceals that 
we have tensored with $\O_{M'}$ on the {\sl left}, for it does
not matter on what side we tensor). Then $\Lambda'$ is
naturally a split almost polynomial algebra of differential
operators on $M'$ (see section 5.1), and the category of
$\Lambda$ modules on $V$ which are schematically supported on
$M'$ is naturally equivalent to the category of
$\Lambda'$-modules on $M'$. This equivalence follows from the
fact that $\Lambda$ necessarily preserves the ideal sheaf of
$M'$ in $V$.  

Now we come back to our given set up, where we apply the above
remark with $N_i$ (or $N_{i,i+1}$) as $V$, $F_i$ (or $F_{i,i+1}$) as
$M$, and the zero section $X_i$ of $N_i$ (or the zero section
$Y^*_i$ of $N_{i,i+1}$) as $M'$.
As usual, let $\D_{N_i}[\log F_i]$ denote the subring of
$\D_{N_i}$ consisting of all operators which preserve the ideal
sheaf of $F_i$ in $\O_{N_i}$. Then we define   
$\D_i$ to be the restriction (see the above remark) of
$\D_{N_i}[\log F_i]$ to the zero section $X_i\subset N_i$, which
is canonically isomorphic to $f_i^*\D_X[\log Y]$. 
Similarly, we define $\D^*_i$ to be the restriction of
$\D_{N_{i,i+1}}[\log F_{i,i+1}]$ to the zero section $Y^*_i$. 

If $H\subset M$ is a nonsingular irreducible component of a
normal crossing divisor $M$ in a nosingular variety $V$, then
the {\bf Euler operator along $H$ in $\D_V[\log M]$} 
is the element $\theta_H \in (\D_V[\log M]|H)$ which 
has the usual definition: if $V$ has local analytic coordinates 
$(x_1,\ldots x_d)$ with $M$ locally defined by $x_1\cdots
x_{d-m}=0$ and $H$ by $x_1=0$, then the operator $\theta_H$ is given 
by the action of $x_1(\partial/\partial x_1)$. 

For each $i\le d-1$ we define the
section $\theta_i$ of $\D^*_i$ as the restriction to $Y^*_i$ of
the Euler operator along $H_{i,i+1}$ in 
$\D_{N_{i,i+1}}[\log F_{i,i+1}]$.  

The above definitions work equally well in the analytic
category. For the remaining basic definitions, we restrict to
the analytic category, with euclidean topology (in particular,
if $T$ was earlier a finite type, reduced scheme over $\C$ then
now the same notation $T$ will denote the corresponding analytic
space). Vector bundles will denote their respective total analytic
spaces. 

Let $U_i$ be the open subset $U_i =N_i-F_i$ of $N_i$, and let $R_i$
be the open subset $N_{i,i+1}-F_{i,i+1}$ of $N_{i,i+1}$. 
Let $N_{i,i+2}$ be the normal bundle to $Z_i$ in $X$, and let
$W_i$ be the open subset of $N_{i,i+2}$ which is the complement 
in $N_{i,i+2}$ of vectors tangent to branches of $Y$.
Similarly, let $N_{i,i+1,i+2}$ be the normal bundle to $Z^*_i$
in $X$. and let $W^*_i$ be the open subset of $N_{i,i+1,i+2}$
which is the complement in $N_{i,i+1,i+2}$ of vectors tangent to
branches of $Y$.

Finally, for $i\le d-1$ we define some central elements $\t_i(c)$ of
certain fundamental groups, which are the topological
counterparts of the operators $\theta_i$. Let $Y^*_i(c)$, where
$c$ varies over $\pi_0(Y^*_i)$, be the
connected components of $Y^*_i$, and let $N_{i,i+1}(c)$,
$F_{i,i+1}(c)$, $H_{i,i+1}(c)$, and $R_i(c)$ be the restrictions
of the corresponding objects to $Y^*_i(c)$.
Let $\tau_i(c)$ be the element in the center of $\pi_1(R_i(c))$,
(with respect to any base point) which is represented by a
positive loop around $H_{i,i+1}(c)$. The fact that 
$\tau_i(c)$ is central, and is unambigously defined, follows
from the following lemma in the topological category.

\begin{lemma}
Let $S$ be a connected topological manifold, and $p:N\to S$ a
complex vector 
bundle on $S$ of rank $r$. Let $F\subset N$ be a closed subset
such that locally over $S^1$, the subset $F$ is the union of $r$
vector subbundles of $N$, each of rank $r-1$, in general
position. Let $H\subset N$ be a vector subbundle of rank $r-1$
such that $H\subset F$. Let $U=N-F$ with projection $p:U\to S$,
which is a locally trivial fibration with fiber $(\C\,^*)\,^r$.
The fundamental group of any fiber is $\Z^r$, with a basis given
by positive loops around the various hyperplanes. 
Let $u_0\in U$ be a base point, and let $p(u_0)=s_0\in S$. Let
$\t_H\in \pi_1(U, u_0)$ be represented by the positive loop around
$H\cap p^{-1}(s_0)$ in the fiber $U_{s_0}$. Then we have the
following: 

(1) The element $\t_H$ is central in $\pi_1(U,u_0)$.

(2) Let $u_1\in U$ be another base point, and let $\t'_H\in
\pi_1(U,u_1)$ be similarly defined. Let $\sigma:[0,1]\to U$ be a
path joining $u_0$ to $u_1$ , and let $\sigma^*: \pi_1(U,u_1)\to
\pi_1(U,u_0)$ be the resulting isomorphism. Then
$\sigma^*(\t'_H)=\t_H$. 
\end{lemma}  

\proof (Sketch) Let $S^1$ be the unit circle with base point $1$,
and let $\gamma: S^1\to U : 1\mapsto u_0$ be another
loop in $U$, based at $u_0$. By pulling back the bundle $N$
under the base chang $p\circ \gamma :S^1\to S$, we can reduce
the statement (1) to the case that $S=S^1$. If base is $S^1$,
then as all complex vector bundles become trivial, the space $U$ 
becomes a product $\C\,^*\times U'$ for some $U'$, and $\t_H$ is
the positive generator of the fundamental group of $\C\,^*$, so is
central in this product. Similarly, (2) follows from a base
change to the unit interval $[0,1]$.

\subsection{Summary of basic notation}

\begin{tabular}{ll}
$X$         & $=$ a nonsingular variety. \\
$d$         & $=$ the dimension of $X$.  \\
$Y$         & $=$ a divisor with normal crossing in $X$.\\
$S_d$       & $=$ $X$ \\
$S_{d-1}$   & $=$ $Y$ \\
            & \\
            & By decreasing induction we define starting with $i=d-2$,\\
$S_i$       & $=$ the singular locus of $S_{i+1}$ for $i\le d-2$. \\
$m$         & $=$ the smallest $i$ for which $S_i$ is nonempty. \\
$I$         & $=$ any nonempty subset of $\{\, m,\ldots ,d-1\,\}$\\
$m(I)$      & $=$ the smallest element of $I$. \\
$p_I:X_I\to S_{m(I)}$       
            & $=$ the finite scheme over $S_{m(I)}$ whose fiber over \\
            & point of $S_{m(I)}$ consists of all nested sequences of 
			  branches \\
            & of $S_j$ at that point, for $j$ varying over $I$. \\
            & In particular when $I=\{\, i\,\}$, we have \\  
$p_i:X_i\to S_i$       
            & $=$ the normalization of $S_i$. \\
            & For nonempty subsets 
			  $J\subset I\subset \{\, m,\ldots ,d-1\,\}$, \\
$X_{i,a}$   & $=$ connected components of $X_i$, \\
            & as $a$ varies over the indexing set $\pi_0(X_i)$ \\ 
$p_{J,I}:X_I\to X_J$
            & $=$ the canonical map. (All these commute.) \\
\end{tabular}

\bigskip
\bigskip

\hfill {\footnotesize (Continued on next page) 

}

\vfill
\pagebreak

\begin{tabular}{ll}			
            & For $m\le i\le d-1$ we put \\
$Y^*_i$     & $=$ $X_{i,i+1}$, the $d-i$ sheeted finite \'etale
              cover of $X_i$ \\ 
            & which splits the branches of $S_{i+1}$ which meet
              along $S_i$ \\ 
   			& For $m\le i\le d-2$ we put \\
$Z_i$       & $=$ $X_{i,i+2}$, the $C^{d-i}_2$ sheeted finite
              \'etale cover of $X_i$ \\ 
            & which splits the branches of $S_{i+2}$ which meet
              along $S_i$ \\ 
$Z^*_i$     & $=$ $X_{i,i+1,i+2}$, the $2$ sheeted finite
              \'etale cover of $Z_i$ \\ 
$f_i:X_i\to X$ 
            & $=$ the composite $X_i\stackrel{p_i}{\to}S_i\to X$ \\

$N_i$       & $=$ the normal bundle to $X_i$ in $X$ under $f_i$,
              defined by \\
			& the exact sequence $0\to T_{X_i}\to f_i^*T_X\to
              N_i\to 0$\\
$T(Y)\subset T_X$
            & $=$ the closed subset of $T_X$ consisting of \\
			& vectors tangent to branches of $Y$    \\
$F_i$       & $=$ the closed subset of $N_i$ which is the image
              of $f^*T(Y)$ \\
            & under the morphism $f^*T_X\to N_i$ of geometric
              vector bundles. \\
$U_i$       & $=$ $N_i-F_i$ \\

$\D_i$      & $=$ the restriction of $\D_{N_i}[\log F_i]$ to the
              zero section $X_i\subset N_i$,\\ 
            & which is canonically isomorphic to $f_i^*\D_X[\log Y]$ \\
$Y^*$       & $=$ $Y^*_{d-1}$ $=$ $X_{d-1}$, the normalization of $Y$ \\ 
$N_{i,i+1}$ & $=$ the normal bundle to $Y^*_i=X_{i,i+1}$ in $X$, \\
            & defined with respect to the morphism $Y^*_i\to
              X_i\to X$ \\
$F_{i,i+1}$ & $=$ the normal crossing divisor in the total space of
              $N_{i,i+1}$, \\
			& defined by vectors tangent to branches of $Y$. \\  
$h_i:Y^*_i\to Y^*$
            & $=$ the canonical map, sending $(x_i,x_{i+1})$ to
              the branch \\
			& $y_{d-1}$ of $Y^*$  which intersects the given
              branch $x_{i+1}$ of $S_{i+1}$ \\ 
			& along given the branch $x_i$ of $S_i$. \\ 			  
$H_{i,i+1}$	& $=$ the hypersubbundle of $N_{i,i+1}$ contained in
              $F_{i,i+1}$, \\
			& defined by vectors tangent to the branch of $Y$ \\
            & given by  $h_i:Y^*_i\to Y^*$ \\		  
$\D^*_i$    & $=$ the restriction of $\D_{N_{i,i+1}}[\log F_{i,i+1}]$ \\
            & to the zero section $Y^*_i\subset N_{i,i+1}$ \\			  
$\theta_i$  & $=$ the Euler operator along $H_{i,i+1}$ in $\D^*_i$. \\
$R_i$       & $=$ the open subset of $N_{i,i+1}-F_{i,i+1}$ of
              $N_{i,i+1}$ which is the \\
			& complement in $N_{i,i+1}$ of vectors tangent to
              branches of $Y$\\ 
$R_i(c)$    & $=$ connected components of $R_i$ as $c$ varies
              over $\pi_0(Y^*_i)$.\\
$\tau_i(c)$ & $=$ the element in the center of $\pi_1(R_i(c))$, which is \\
            & given by a positive loop around $H_{i,i+1}(c)$\\
$N_{i,i+2}$ & $=$ the normal bundle to $Z_i$ in $X$. \\
$W_i$       & $=$ the open subset of $N_{i,i+2}$ which is the complement \\
            & in $N_{i,i+1,i+2}$ of vectors tangent to branches of $Y$\\
$N_{i,i+1,i+2}$ 
            & $=$ the normal bundle to $Z^*_i$ in $X$. \\
$W^*_i$     & $=$ the open subset of $N_{i,i+1,i+2}$ which is the
              complement \\ 
            & in $N_{i,i+1,i+2}$ of vectors tangent to branches of $Y$\\
\end{tabular}

\vfill

\section{Pre-$\D$-modules on $(X,Y)$} 
In this section, we first define the notion of a
pre-$\D$-module. Then we consider the special case when $X$ is a
polydisk. Finally, we give some historical motivation.

\subsection{Global definition}
Let $X$ is any smooth variety and $Y$ a normal 
crossing divisor. We follow the notation introduced in section 2.
The following definitnion works equally well in the algebraic or
the analytic categories.

\begin{definition}\rm 
A {\bf pre-$\D$-module} $\EE = (E_i,t_i,s_i)$ on 
$(X,Y)$ consists of the following. 

(1) For each $m\le i\le d$, $E_i$ is a vector bundle on $X_i$ 
(of not necessarily constant rank) together with a 
structure of $\D_i$-module, 

(Note that by (1), for each $m\le i\le d-1$, the pullbacks 
$E_{i+1}| Y^*_i$ and $E_i| Y^*_i$ under the respective maps 
$Y^*_i\to X_{i+1}$ and $Y^*_i\to X_i$ have a natural structure
of a $\D^*_i$-module.)

(2)  For each $m\le i\le d-1$, 
$t_i:(E_{i+1}| Y^*_i) \to (E_i| Y^*_i)$ and 
$s_i:(E_i |Y^*_i) \to (E_{i+1} | Y^*_i)$ are 
$\D^*_i$-linear maps, such that 

\begin{eqnarray*}
s_it_i &=& \theta_i {\mbox{~{\rm on}~}} 
                     E_{i+1}|Y^*_i \\ 
t_is_i &=& \theta_i {\mbox{~{\rm on}~}} 
                     E_i| Y^*_i \\ 
\end{eqnarray*}

(3) Let $m\le i \le d-2$. 
Let $\pi:Z^*_i\to Z_i$ be the projection
$$p_{\{\,i,i+2\,\} ,\{\,i,i+1,i+2\,\} }: 
X_{\{\,i,i+1,i+2\,\} }\to X_{\{\,i,i+2\,\} }$$
Let $E_{i+2}|Z_i$ and $E_i|Z_i$ be the pullbacks of $E_{i+2}$ and $E_i$ 
under respectively the composites
$Z_i=X_{\{\,i,i+2\,\} }\to X_{i+2}$ and 
$Z_i=X_{\{\,i,i+2\,\} }\to X_i$. 
{\rm (Note that there is no object called $E_{i+1}|Z_i$.)} 
We will denote 
the pullback of $E_{i+1}$ under 
$Z^*_i=X_{\{\,i,i+1,i+2\,\} } \to X_{i+1}$ by $E_{i+1}|Z^*_i$.
We will denote $\pi^*(E_{i+2}|Z_i)$ by $E_{i+2}|Z^*_i$ and 
$\pi^*(E_i|Z_i)$ by $E_i|Z^*_i$.
Let 
$$a_{i+2}:E_{i+2}|Z_i \to \pi_*\pi^*(E_{i+2}|Z_i) = \pi_*(E_{i+2}|Z^*_i)$$
$$a_i:E_i|Z_i \to \pi_*\pi^*(E_i|Z_i) = \pi_*(E_i|Z^*_i)$$
be adjunction maps, and let the cokernels of these maps be denoted by 
$$q_{i+2}:\pi_*(E_{i+2}|Z_i)\to Q_{i+2}$$
$$q_i:\pi_*(E_i|Z_i)\to Q_i$$
Then we impose the
requirement that the composite map 
$$E_{i+2}|Z_i\stackrel{a_{i+2}}{\to} \pi_*(E_{i+2}|Z^*_i) 
\stackrel{\pi_*(t_{i+1}|Z^*_i)}{\to}  
\pi_*(E_{i+1}\vert Z^*_i) 
\stackrel{\pi_*(t_i|Z^*_i)}{\to} \pi_*(E_i\vert Z^*_i)
\stackrel{q_i}{\to} Q_i  $$
is zero.

(4) Similarly, we demand that for all $m\le i\le d-2$ the
composite map 
$$Q_{i+2}\stackrel{q_{i+2}}{\leftarrow} \pi_*(E_{i+2}\vert Z^*_i) 
\stackrel{\pi_*(s_{i+1}|Z^*_i)}{\leftarrow} 
\pi_*(E_{i+1}\vert Z^*_i) 
\stackrel{\pi_*(s_i|Z^*_i)}{\leftarrow} \pi_*(E_i\vert Z^*_i) 
\stackrel{a_i}{\leftarrow} E_i|Z_i$$
is zero.

(5) Note that as $\pi:Z^*_i\to Z_i$ is a 2-sheeted cover, 
for any sheaf $\F$ on $Z_i$ the new sheaf $\pi_*\pi^*(\F)$ on $Z_i$ has 
a canonical involution coming from the deck transformation for 
$Z^*_i\to Z_i$ which transposes the two points over any base point. 
In particular, the
bundles $\pi_*(E_{i+2}|Z^*_i)=\pi_*\pi^*(E_{i+2}|Z_i)$ and 
$\pi_*(E_i|Z^*_i)=\pi_*\pi^*(E_i|Z_i))$ have canonical
involutions, which we denote by $\nu$. We demand that the
following diagram should commute. 

Diagram III.
$$\begin{array}{ccccc}
\pi_*(E_{i+1} | Z^*_i) & \stackrel{\pi_*(s_{i+1}|Z^*_i)}{\to} & 
\pi_*(E_{i+2}|Z^*_i) & \stackrel{\nu}{\to} & \pi_*(E_{i+2}|Z^*_i) \\
{\scriptstyle \pi_*(t_i|Z^*_i)}\downarrow & & & & 
\downarrow {\scriptstyle \pi_*(t_i|Z^*_i)}\\ 
\pi_*(E_i|Z^*_i) & \stackrel{\nu}{\to} & \pi_*(E_i|Z^*_i)  
& \stackrel{\pi_*(s_i|Z^*_i)}{\to} & \pi_*(E_{i+1} | Z^*) \\
\end{array}$$

A {\bf homomorphism} $\varphi :\EE \to \EE'$ of pre-$\D$-modules
consists of a collection $\varphi_i:E_i\to E'_i$ of
$\D_i$-linear homomorphisms which make the obvious diagrams commute. 
\end{definition}

\rem As the adjunction maps are injective (in particular as
$a_i$ is injective),  
the condition (3) is equivalent to demanding the existence of a
unique $f$ which makes  the following diagram commute.

Diagram I.
$$\begin{array}{ccccc}
E_{i+2}|Z_i & & \stackrel{f}{\longrightarrow} & & E_i|Z_i \\
a_{i+2}\downarrow & & & & \downarrow a_i\\
\pi_*(E_{i+2}\vert Z^*_i) &
\stackrel{\pi_*(t_{i+1}|Z^*_i)}{\to} & 
\pi_*(E_{i+1}\vert Z^*_i) &  
\stackrel{\pi_*(t_i|Z^*_i)}{\to} & 
\pi_*(E_i\vert Z^*_i)\\
\end{array}$$
Similarly, the condition (4) is equivalent to the existence of a
unique homomorphism $g$ which makes 
the following diagram commute.

Diagram II.
$$\begin{array}{ccccc}
E_{i+2}|Z_i & & \stackrel{g}{\longleftarrow} & & E_i|Z_i \\
a_{i+2}\downarrow & & & & \downarrow a_i\\
\pi_*(E_{i+2}\vert Z^*_i) &
\stackrel{\pi_*(s_{i+1}|Z^*_i)}{\leftarrow} & 
\pi_*(E_{i+1}\vert Z^*_i) &  
\stackrel{\pi_*(s_i|Z^*_i)}{\leftarrow} & 
\pi_*(E_i\vert Z^*_i)\\
\end{array}$$

\rem\label{local} The above definition is local, in the sense that 

(1) if $U_{\alpha}$ is an open covering of $X$ in the algebraic
or analytic category, and $\EE_{\alpha}$ is a collection of
pre-$\D$-modules on $(U_{\alpha},\,Y\cap U_{\alpha})$ together
with isomorphisms
$\varphi_{\alpha,\beta}:(\EE_{\beta}|U_{\alpha, \beta})
\to (\EE_{\alpha}|U_{\alpha, \beta})$ which form a
$1$-cocycle, then there exists a unique (upto unique isomorphism) 
pre-$\D$-module $\EE$ on $(X,Y)$ obtained by gluing.

(2) If $\EE$ and $\EE'$ are two pre-$\D$-modules on $(X,Y)$,
then a homomorphism $\varphi :\EE\to \EE'$ is uniquely defined
by a collection of homomorphisms over $U_{\alpha}$ which match
in $U_{\alpha, \beta}$.

\begin{definition}\label{goodresieigen}\rm
Consider the action of $\theta_i$ on $E_{i+1}|Y^*_i$, which is
$\O_{Y^*_i}$-linear. Note that in the global algebraic case,
compactness of $Y^*_i$ implies that 
the characteristic polynomial of the resulting endomorphism 
of $E_{i+1}|Y^*_i$ is constant on each component of $Y^*_i$. 
We say that a {\bf pre-$\D$-module $\EE=(E_i,s_i,t_i)$ has good
residual eigenvalues} if for each $i \le d-1$ no two eigenvalues of
$\theta_i\in End(E_{i+1}|Y^*_i)$, on any two components of
$Y^*_i$ which map down to intersecting subsets of $X$, 
differ by a non-zero integer.  
\end{definition}

\rem Note that the above definition does not prohibit two
eigenvalues of $\theta_i$ on $E_i|Y^*_i$ from differing by
non-zero integers. Also, note that the definition can involve more
than one component of $Y^*_i$ at a time: it is stronger than
requiring that on each component of $Y^*_i$ no two eigenvalues
should differ by nonzero integes.

\subsection{Restriction to a polydisk}

There exists an open covering of $X$, where each open subset is
a polydisk in $\C^d$ with coordinates $x_i$, defined by $|x_i|<
1$, whose intersection with $Y$ is defined by  
$\prod_{i\le r}x_i=0$ for a variable integer $r\le d-m$. 
It is possible globally that the irreducible components of $Y$
are singular. Moreover, it is possible that various branches of 
$Y$ meeting at a point get interchanged as one moves around. 
This does not happen in a polydisk of the above kind, so the
definition of a pre-$\D$-module becomes much simpler. We give it
in detail in view of remark \ref{local} above.

Let $X$ be a polydisk in $\C^d$ around the origin, with coordinates 
$x_1,\ldots, x_d$, and 
let $m$ be some fixed integer with $0\le m\le d$.
If $m\le d-1$, let $Y\subset X$ be the normal crossing divisor 
defined by $\prod _{1\le i\le d-m}x_i=0$. 
If $m=d$, we take $Y$ to be empty.

We will follow the notation summarised in section 2. 
We have a filtration $Y=S_{d-1}\supset \ldots \supset S_m$ 
where each $S_i$ is the singularity set of $S_{i+1}$ for $m\le
i\le d-2$, and $S_m$ is nonsingular. 
Note that the irreducible components 
of $S_i$ are as follows. For any subset $A\subset \{\,1,\ldots,d-m\,\}$ of
cardinality $d-i$, we have a component $S_A$ of $S_i$ defined by the ideal
generated by all $x_j$ for $j\in A$ (total $C^{d-m}_{d-i}$ components). 
Then the normalization $X_i$ of $S_i$ is simply the disjoint union of all
the $S_A$. Therefore $X_{i,A} = S_A$, which are polydisks of dimension 
$d- |A|$, are the components of $X_i$. Whenever $k\in A$, we
have an inclusion $X_A\to X_{A-\{ k\}}$. This is identified in
our earlier notation with a component of $Y^*_{i-1} \to X_i$.

It follows from its general definition 
that a pre-$\D$-module $\EE$ on $(X,Y)$ 
consists of the following data. 

(1) For each $A\subset \{\, 1,\ldots , d-m\,\}$, we are given a
vector bundle $E_A$ on $S_A$, together with the structure of a
$\D_X[\log Y]$-module. The $E_A$ with $|A|=d-i$ are the
restrictions of $E_i$ to the components $S_A$ of $X_i$.

(2) For any $k \in A\subset \{\, 1,\ldots , d-m\,\}$, we have
$\D_X[\log Y]$-linear homomorphisms  
$t^k_A:E_{A-k}|X_A \to E_A$ and $s^k_A:E_A \to E_{A-k}|X_A$ 
such that 
$$s^k_At^k_A = x_k\partial/\partial x_k {\mbox{~{\rm on}~}} E_{A-k}|X_A$$
$$t^k_As^k_A =  x_k\partial/\partial x_k {\mbox{~{\rm on}~}} E_A$$
In terms of earlier notation, the $t^k_A$ (respectively, the
the $s^k_A$) with $|A|=d-i$ make up the $t_i$ (respectively the
$s_i$).

(3) Let $k\ne \ell$ such that $k,\, \ell \in A\subset \{\,
1,\ldots , d-m\,\}$. Then we must have 

\begin{eqnarray*}
t^k_At^{\ell}_{A-k} = t^{\ell}_At^k_{A-{\ell}}\\
s^k_{A-{\ell}}s^{\ell}_A = s^{\ell}_{A-k}s^k_A \\
t^k_{A-{\ell}}s^{\ell}_{A-k} = s^{\ell}_A t^k_A \\
\end{eqnarray*}

The above three equations respectively embody the conditions
that diagrams I, II, and III in the definition of a
pre-$\D$-module must commute.

\subsection{Motivation for the definition}

The definition of a pre-$\D$-module may be regarded as
another step in the programme of giving concrete representations
of regular holonomic $\D$ modules and perverse sheaves. The
earlier steps relevant to us are the following.

{\bf (1)}  Deligne's description (1982) of a perverse sheaf on a
disk with singularity at the origin, in terms of pairs of vector spaces
and linear maps and Malgrange's description of corresponding
regular holonomic $\D$-modules.  

{\bf (2)} Verdier's functor of specialization (Asterisque
101-102), and his description of extension of a perverse sheaf
across a closed subspace (Asterisque 130).  

{\bf (3)} Similar construction by Malgrange for regular
holonomic $\D$-modules in place of perverse sheaves.

{\bf (4)} Verdier's description of a perverse sheaf on the total
space of a line bundle $L$ on a smooth variety $S$, in terms of
two local systems on $L-S$ ($=$ the com\-ple\-ment of the zero
section) and 
maps between them (Asterisque 130, 1985). 

{\bf (5)} Definition of a pre-$\D$-module on $(X,Y)$ when $Y$ is
nonsingular, which can be obtained by choosing compatible
logarithmic lattices in a combination of Step 3 and Step 4  
(see [N-S]).

{\bf (6)} Description by Galligo, Granger, Maisonobe of perverse
sheaves on a polydisk with coordinates $(z_1,\ldots,z_n)$ with
respect to the smoothening stratification induced by the normal
crossing divisor $z_1\cdots z_n=0$, in terms of a hypercube of
vector spaces and linear maps (1985). 

The description (6) for a
polydisk with a normal cross\-ing divisor is local and coordinate
dependent like the description (1) for the disk with a point. One first 
makes it coordinate free and globalizes it in order to have the equivalent
of (4) (which gives us finite descriptions of perverse sheaves
described in section 6 below), and then puts level structures
generalizing (5) to arrive at the above definition of a pre-$\D$-module.

One of the problems in globalizing the local hypercube description
is that one can not unambiguously label the branches of $Y$ which meet
at a point, because of twistedness of the divisor. This is
taken care of by normalizing the closed strata $S_i$ and 
going to the coverings $\pi :Z^*_i\to Z_i$. 

The requirement that the various composites of $s$ and $t$ should
give endomorphisms expressible in terms of Euler vector fields
is present in the local hypercube description in much the same form.

The three commutative diagrams I, II, and III in the definition
respectively embody the globalizations of the 
conditions in the hypercube description that 

I: the two canonical maps $C_x$ and $C_y$ should commute,

II: the two variation maps $V_x$ and $V_y$ should commute, and 

III: we should have $C_xV_y=V_yC_x$.

\bigskip

\section{From pre-$\D$-modules to $\D$-modules}
In section 4.1, we directly describe the $\D$-module
associated to a pre-$\D$-module in the special case where $Y$ is
nonsingular. In section 4.2, we will associate a $\D$-module
to a pre-$\D$-module when $Y$ is normal crossing. This is done
by first doing it on polydisks, and then patching up. Finallly,
in section 4.3 we show how to find a pre-$\D$-module with
good residual eigenvalues over a given $\D$-module, proving that
the functor from pre-$\D$-modules on $(X,Y)$ with good residual
eigenvalues to regular holonomic $\D$-modules on $X$ whose
caracteristic variety is contained in $N^*(Y)$ is essentially
surjective.  

\subsection{The case when $Y$ is smooth}

First we treat the case where the divisor $Y$ is nonsingular.
To a pre-$\D$-module $(E,F,t,s)$ on $(X,Y)$, where
$E=E_d$ is a logarithmic connection on $(X,Y)$, $F=E_{d-1}$ is a vector 
bundle on $Y$ with structure of a $\D_X[\log Y]$-module, and 
$t:E|Y\to F$ and $s:F\to E|Y$ are $\D_X[\log Y]$-linear maps with 
$st = \theta_Y$ on $E|Y$, and $ts=\theta_Y$ on $F$, we will directly 
associate the following $\D$-module $M$ on $X$. (This was 
indirectly described in [N-S]).

Let $M_0=E$, and let $E\oplus_sF$ denote the subsheaf of $E\oplus F$
consisting of sections $(e,f)$ such that $(e|Y)=s(f)$. Let $\O_X(Y)$ be 
the line bundle on $X$ defined by the divisor $Y$ as usual, and let
$M_1=\O_X(Y)\otimes (E\oplus_sF)$. Let $M_0\hookrightarrow M_1$ be the 
inclusion defined by sending a local section $e$ of $M_0=E$ to the local 
section $(1/x)\otimes (xe,0)$, where $x$ is a local generator for the 
ideal of $Y$ in $X$ (this can be readily seen to be independent of the 
choice of $x$).

We make $E\oplus_sF$ is a $\D_X[\log Y]$-module by putting for any 
local section $\xi$ of $T_X[\log Y]$ and $(e,f)$ of $E\oplus_sF$,
$$\xi (e,f) = (\xi(e) , \, t(e|Y) + \xi(f))$$
The right hand side may again be checked to be in $E\oplus_sF$, using the 
relation $st=\theta_Y$ on $E|Y$. As $\O_X(Y)$ is naturally a $\D_X[\log 
Y]$-module, this now gives the structure of a (left) $\D_X[\log 
Y]$-module on the tensor product $M_1=\O_X(Y)\otimes_{\O_X}(E\oplus_sF)$. 
Moreover, the inclusion $M_0\hookrightarrow M_1$ defined above is 
$\D_X[\log Y]$-linear. 

We now define a connection $\nabla : M_0 \to \Omega^1_X\otimes
M_1$ by putting, for any  
local sections $\eta$ of $T_X$ and $e$ of $M_0$,
$$\eta (e) = (1/x)\otimes ((x\eta)(e),\, \eta(x) t(e|Y))$$
where $x$ is any local generator of the ideal of $Y$.
The right hand side makes sense because $x\eta$ is a section of $T_X[\log 
Y]$, and so $(x\eta)(e)$ is defined by the logarithmic connection on $E$. 
It can be checked that the above formula is independent of the choice of 
$x$, is $\O_X$-linear in the variable $\eta$, and the resulting map 
$\nabla: M_0 \to \Omega^1_X\otimes M_1$ satisfies the Leibniz rule.
Moreover, the following diagram commutes, where the maps $M_i\to 
\Omega^1_X[\log Y]\otimes M_i$ (for $i=0$ and for $i=1$) are given by the 
$\D_X[\log Y]$-module structure on $M_i$.

$$\begin{array}{ccc}
M_0        & \to & \Omega^1_X\otimes M_1 \\
\downarrow &     & \downarrow \\
M_1        & \to & \Omega^1_X[\log Y]\otimes M_1  \\
\end{array}$$

Now let $M$ be the $\D_X$-module which is the quotient of 
$\D_X\otimes_{\D_X[\log Y]}M_1$ by the submodule generated by elements of 
the type $\eta\otimes e - 1\otimes \eta(e)$ where $e$ is a local section 
of $E=M_0\subset M_1$ and $\eta$ is a local section of $T_X$. Then $M$ is 
the $\D_X$-module that we associate to the pre-$\D$-module $(E,F,t,s)$.

\centerline{\sl Relation with $V$-filtration}

We now assume that the generalized eigenvalues of 
$\theta_Y$ on $E|Y$ do not differ by nonzero integers  
(it is actually enough to assume this along each connected
component of $Y$,  
but for simplicity we will assume that $Y$ is connected). Let $\mu$ 
the only possible integral eigenvalue (when there are more
components in $Y$, 
there can be a possibly different $\mu$ along each component). 
Under this assumption,
We now construct a $V$-filtration on $M$ along the divisor $Y$. Put
$V^{\mu}(M)$ to be the image of $M_0$ and $V^{\mu +1}(M)$ to be the 
image of $M_1$ in $M$. For $k\ge 1$ put 
$V^{\mu -k}(M)$ to be the image of $I_Y^kM_0 \subset M_0$ and
$V^{\mu +k}(M)$   
to be the image of $\O_X((k-1)Y)\otimes M_1$ in $M$. Then by
definition each  
$V^kM$ for $k\in \Z$ is an $\O_X$-coherent $\D_X[\log Y]$-module, with 
$\eta (V^k(M))\subset V^{k+1}(M)$ and $I_YV^{k+1}(M) \subset V^k(M)$ 
for all $k$. 

Conversely, let $M$ be a regular holonomic $\D$-module on $X$ with 
$car(M)\subset N^*(Y)$. Let $V^k(M)$ be a $V$-filtration 
with $\mu$ the only integer in the fundamental domain chosen for the 
exponential map $\C \to \C\,^*$). Then put $E=V^{\mu}(M)$, 
$F = N^*_{Y,X}\otimes (V^{{\mu}+1}(M)/V^{\mu}(M))$, 
$t:(E|Y)\to F$ is defineded by putting
$$t(e|Y) = x\otimes (\partial/\partial x)e$$ 
where $x$ is a local 
generator of $I_Y$ (which is independent of the choice of $x$), 
and $s:F\to E|Y$ defined by simply the multiplication 
$I_Y\times V^{\mu +1}M\to V^{\mu}M$ (note for this that
$N^*_{Y,X}=I_Y/I_Y^2$). Then we get a pre-$\D$-module
$(E,F,t,s)$. Given a choice of a fundamental domain for the
exponential map, the above two processes are inverses of each other.

\subsection{General case of a normal crossing $Y$}

Let $\EE =(E_i,t_i,s_i)$ be a pre-$\D$-module on $(X,Y)$, 
where we now allow $Y$ to 
have normal crossings. Let the sheaf $F$ on $X$ be the subsheaf of 
$\oplus (p_i)_*E_i$ whose local sections consist of all tuples
$(e_i)$ where  
$e_i\in (p_i)_*E_i$, such that $s_i(e_i|Y^*_i)=e_{i+1}|Y^*_i$. 
This is a $\D_X[\log Y]$-submodule of $\oplus (p_i)_*E_i$ as may be seen. 
Let $G=\O_X(Y)\otimes F$. As $\O_X(Y)$ is naturally a $\D_X[\log
Y]$-module, $G$ has a natural structure of a left $\D_X[\log Y]$-module.
The $\D_X$-module $M$ that we are going to associate to the pre-$\D$-module 
$\EE$ is going to be a particular quotient of the $\D_X$-module 
$\D_X\otimes_{\D_X[\log Y]}G$. 

\medskip

\centerline{\sl The case of a polydisk}

Let $x_1,\ldots,x_d$ be coordinates on the polydisk, let $0\le r\le d$
and let $Y$ be defined by the polynomial $P(x)=\prod_{1\le k\le
r}x_k$ (in particular $P=1$ and therefore $Y$ is 
empty if $r=0$). A pre-$\D$-module $\EE$ on $X$ has been described 
already in section 3.1 above. 

For each $k$ such that $1\le k\le r$, we define a subsheaf 
$F_k\subset F$ (where $F$ is the submodule of $\oplus
(p_A)_*E_A$ defined above using the $s_A$) as follows. 
$$F_k = x_kF\subset F$$
This defines an $\O_X$-coherent module, which is in fact a $\D_X[\log 
Y]$-submodule of $F$ as may be checked. We now put $G_k=\O_X(Y)\otimes 
F_k \subset \O_X(Y)\otimes F = G$. This is therefore a $\D_X[\log 
Y]$-submodule of $G$. We now define the operator
$$\partial/\partial x_k:G_k \to G$$ as follows.
If $e=(1/P)\otimes (e_A)$ is a section of $G_k$, we put 
$\partial_k(e) = (1/P)\otimes (f_B)$ where 
$$f_B = (x_k\partial/\partial x_k)(e_B) + t^k_B(e_{B-k})$$
where by convention $t^k_B=0$ whenever $k$ does not belong to $B$.

Let $K$ be the $\D_X$-submodule of $\D_X\otimes_{\D_X[\log Y]}G$ 
generated by elements of the form 
$$\partial_k\otimes e -1\otimes \partial_k(e)$$
where $e\in G_k$. 

\begin{lemma}
The submodule $K$ is independent of the choice of local coordinates, and 
restricts to the corresponding submodule on a smaller polydisk.
\end{lemma}

\proof This is a local coordinate calculation, using the chain rule of 
partial differentiation under a change of coordinates. We omit the 
details. 

\medskip

\centerline{\sl Back to the global case}

By the above lemma applied to an open covering of $X$ by
polydisks, we get a globally defined submodule $K$ of
$\D_X\otimes G$ (where the later is already defined globally).
We now put $M$ to be the quotient of $\D_X\otimes G$ by $K$. This is our 
desired $\D_X$-module.

\subsection{The $V$ filtration for a polydisk}

Let $X$ be a polydisk with coordinates
$x_1,\ldots,x_d$ and $Y$ be defined by $\prod_{k\in\Lambda}x_k$ for
some initial subset ${\Lambda}=\{\,1,\ldots,r\,\}\subset
\{\,1,\ldots,d\,\}$. 

Let $M$ be a regular holonomic $\D_X$-module on $(X,Y)$.
We now describe a pre-$\D$-module $\EE$ such that $M$ is
associated to it. This is done via a $V$-filtration of $M$.
We assume that we have fixed some fundamental domain 
$\Sigma$ for $exp:\C\to \C\,^*:z\mapsto e^{2\pi iz}$, in order to
define the $V$-filtration. Let $\mu\in \Sigma$ be the only integer. 

Note that the filtration will be by sub $\D_X[\log Y]$-modules 
$V^{(n_1,\ldots,n_r)}$ where $n_k\in \Z$, partially ordered as follows.
If $(n_1,\ldots,n_r)\le (m_1,\ldots,m_r)$ 
(which means $n_k\le m_k$ for each $k\in {\Lambda}$), 
then $V^{(n_1,\ldots,n_r)} \subset V^{(m_1,\ldots,m_r)}$. 
In particular, there is a portion of the filtration 
indexed by the power set of ${\Lambda}$, where for any subset
$A$ of ${\Lambda}$, we put 
$V^A = V^{(n_1,\ldots,n_r)}$ where $n_k=\mu +1$ if $k\in A$ and
$n_k=\mu$ otherwise. 
By this convention, note that $V^{\phi} =V^{(\mu ,\ldots,\mu )}$ for the 
empty set $\phi$, and $V^{\Lambda}=V^{(\mu+1,\ldots,\mu +1)}$.

Note that if $k\in A$ then multiplication by $x_k$ defines a map 
$$x_k:V^AM\to V^{A-k}M$$
and differentiation by $\partial/\partial x_k$ defines a map 
$$\partial_k:V^{A-k}M \to V^AM$$
Let for any nonempty $A$, 
$$gr^AM = {V^AM \over \sum_B V^BM}$$
where $B$ varies over all proper subsets of $A$.
All eigenvalues of $x_k\partial/\partial x_k$ on $gr^AM$ lie in $\Sigma$.

The above defines maps
$$x_k: gr^AM\to gr^{A-k}M$$
and
$$\partial_k : gr^{A-k}M\to gr^AM$$

Then we can define a 
pre-$\D$-module $\EE=(E_A,t^k_A,s^k_A)$ as follows:
The $\D_X[\log Y]$-modules $E_A$ are defined by 
$$E_A = \O_X(-\sum_{k\in A} Y^k)\otimes_{\O_X}gr^AM$$
where $Y^k\subset X$ is the divisor defined by $x_k=0$.
Let $k\in A$. To define $t^k_A$ and $s^k_A$, note that the
differentiation $\partial_k: gr^{A-k}M\to gr^AM$ induces a map 
$\O_X(Y^k)\otimes gr^{A-k}M\to gr^AM$, as $\partial_k$ can be canonically 
identified with the section $1/x_k$ of $\O_X(Y^k)$. 
Tensoring this by the identity map on 
$\O_X(-\sum_{\ell\in A}Y^{\ell})$, and restricting to $X_A$, we get 
$t^k_A: (E_{A-k}|X_A)\to E_A$. Also, the map $x_k: gr^AM\to gr^{A-k}M$ 
defines a homomorphism $\O_X(-Y^k)\otimes gr^AM\to gr^{A-k}M$
which after tensoring by the identity map on $\O_X(-\sum_{\ell
\in A-k}Y^{\ell})$ gives  
$s^k_A: E_A\to (E_{A-k}|X_A)$. 

Then it can be checked that we indeed get a pre-$\D$-module $\EE$,
such that it has good residual eigenvalues, lying in $\Sigma$,
such that $M$ is the $\D$-module associated to it. 
It can again be checked that the above proceedure over a
polydisk is coordinate independent, so glues up to give such a
correspondence globally over $(X,Y)$.

\rem If we began with the $\D$-module $M$ associated to a
pre-$\D$-module $\EE$ with good residual eigenvalues 
lying in $\Sigma$, then the above will give back $\EE$, as then
$V^{\phi}M = E_d$ and for $A\subset \Lambda$ we will get 
$$V^AM = (\prod_{\ell \in {\Lambda}-A}x_{\ell}) G \subset G$$ 
where $G=\O_X(Y)\otimes F$ where $F\subset \oplus_BE_B$ as above.

\rem (The [G-G-M]-hypercube for a polydisk) :
Let $W_A$ be the fiber of $E_A$ at the origin of the polydisk.
For $k\in A$  
let $t^k_A : W_{A-k}\to W_A$ again denote the restriction of 
$t^k_A :E_{A-k} \to E_A$ to the fiber at origin, and let 
$v^k_A: W_A\to W_{A-k}$ be defined by the following formula.
$$v^k_A = {\exp(2\pi i\theta_k)-1\over \theta_k} s^k_A$$
Then $(W_A,t^k_A,v^k_A)$ is the hypercube description of $M$ as
given by  
Galligo, Granger, Maisonobe in [G-G-M].

\bigskip

\section{Moduli for semistable pre-$\D$-modules}

In this section we define the concepts of semistability and
stability for pre-$\D$-modules, and construct a coarse moduli.
The main result is Theorem \ref{maintheorem} below. 

\subsection{Preliminaries on $\Lambda$-modules}

Simpson has introduced a notion of modules over rings of
differential operators which we first recall (see section 2 of [S]).

Let $X$ be a complex scheme, of finite type over $\C$, and let
$\Lambda$ be a sheaf of $\O_X$-algebras (not necessarily
non-commutative),  
together with a filtration by subsheaves of abelian groups
$\Lambda_0\subset \Lambda_1 \subset \ldots \Lambda$ which
satisfies the following properties.  

(1) $\Lambda = \cup \Lambda_i$ and $\Lambda_i\cdot \Lambda_j
\subset \Lambda_{i+j}$. (In particular, $\Lambda_0$ is a
subring, and each $\Lambda_i$ is a $\Lambda_0$-bimodule.)

(2) The image of the homomorphism $\O_X\to \Lambda$ is equal to
$\Lambda_0$. (In particular, each $\Lambda_i$ is an
$\O_X$-bimodule). 

(3) Under the composite map $\C_X\hookrightarrow \O_X \to \Lambda$, the
image of the constant sheaf $\C_X$ is contained in the center of
$\Lambda$.
 
(4) The left and right $\O_X$-module structures on the $i$th
graded piece $Gr_i(\Lambda) = \Lambda_i/\Lambda_{i-1}$ are
equal. 

(5) The sheaves of $\O_X$-modules $Gr_i(\Lambda)$ are coherent.
(In particular, each $\Lambda_i$ is bi-coherent as a bi-module
over $\O_X$, and their union $\Lambda$ is bi-quasicoherent.)

(6) The associated graded $\O_X$-algebra $Gr(\Lambda)$ is
generated (as an algebra) by the piece $Gr_1(\Lambda)$. 

(7) (`Split almost polynomial' condition) : The homomorphism
$\O_X\to \Lambda_0$ is an isomorphism, the $\O_X$-module
$Gr_1(\Lambda)$ is locally free, the graded ring
$Gr(\Lambda)$ is the symmetric algebra over $Gr_1(\Lambda)$, and
we are given a left-$\O_X$-linear splitting $\xi
:Gr_1(\Lambda)\to \Lambda_1$ for the left-$\O_X$-linear
projection $\Lambda_1\to Gr_1(\Lambda)$. 

\rem The condition (7) is not necessary for the moduli construction,
but allows a simple description (see lemma 2.13 of [S]) of the
structure of $\Lambda$-module, just as the structure of a
$\D_X$-module on an $\O_X$-module can be described in terms of the
action of $T_X$. 

The pair $(\Lambda, \xi)$ is called a {\bf split almost
polynomial algebra of differential operators} on $X$. Simpson
also defines this in the relative situation $X\to S$, and treats
basic concepts such as base change, which we will assume.  

A {\bf $\Lambda$-module} will always mean a left $\Lambda$-module
unless otherwise indicated. For any complex scheme $T$, a {\bf
family} $E_T$ of $\Lambda$-modules parametrized by $T$ has an
obvious definition (see [S]). 

The following basic lemma is necessary to parametrize families
of pre-$\D$-modules.  

\begin{lemma}\label{nnreplem}
({\bf Coherence and representability of integrable direct images})
Let $\Lambda$ be an algebra of differential operators on $X$. 
Let $E_T$ and $F_T$ be a families of $\Lambda$-modules on $X$
parametrised by a scheme $T$ (which is locally noetherian and of
finite type over the field of complex numbers).

(i) The sheaf $(\pi _T)_*\underline{Hom}_{\Lambda_T}(E_T,F_T)$ is a
coherent sheaf of $\O_T$ modules. 

(ii) Consider the contravarient functor
from schemes over $T$ to the category of abelian groups, which associates
to ${T'}\lra T$ the  abelian group $Hom_{\Lambda_{T'}}(E_{T'}),F_{T'})$ 
where $E_{T'}$ and $F_{T'}$ are the pullbacks under $X\times {T'}\to
X\times T$. Then there exists a linear scheme $V\lra T$ which
represents this functor.  
\end{lemma}

\proof The above lemma is a stronger version of lemma 2.7 in
[N]. The first step in the proof is the following lemma, which
is an application of the Grothendieck complex of
semi-\-continuity theory. 

\begin{lemma}\label{ega} {\rm (EGA III 7.7.8 and 7.7.9)}
Let $Z\to T$ be a projective morphism where $T$ is noetherian,
and let $\F$ and $\G$ be coherent sheaves on $Z$ such that $\G$
is flat over $T$. Consider the contra functor $\varphi$ from the
category of schemes over $T$ to the category of abelian groups, which
associates to any $T'\to T$ the abelian group of all
$\O_{Z\times_TT'}$-linear homomorphisms from $\F_{T'}$ to
$\G_{T'}$. Then $\varphi$ is representable by a linear scheme
$W$ over $T$. 
\end{lemma}

Now we prove lemma \ref{nnreplem}. Forgetting the structure of
$\Lambda$-modules on $E_T$ and $F_T$ and treating them just as
$\O$-modules, let $W\to T$ be the linear scheme given by the
above lemma \ref{ega}. Then $W$ parametrizes a universal family
of $\O_{X\times W}$-linear morphism $u:E_W\to F_W$. The
condition of $\Lambda_W$-linearity on $u$ defines a closed
linear subscheme $V$ of $W$. 
By its construction, for any base change $T'\to T$, we have
canonical isomorphism $Mor_T(T,V) = Hom_{\Lambda_{T'}}(E_{T'},F_{T'})$ 
which proves the lemma \ref{nnreplem}.

\subsection{Families of pre-$\D$-modules}

We now come back to $(X,Y)$ as before, and our earlier notation.
Recall that $X_i$ is the normalization of $S_i$, which is 
$i$-dimensional, nonsingular if non-empty. Let $X_{i,a}$, as $a$
varies over the indexing set $\pi_0(X_i)$, be the connected
components of $X_i$. We denote by $\D_{i,a}$ the restriction of
$\D_i$ to the component $X_{i,a}$ of $X_i$.

\begin{lemma}\label{lambda}
Each $\D_{i,a}$ satisfies the above properties (1) to (7) (where in
(7), we take $\xi : (\D_{i,a})_1\to \O_{X_{i,a}}$ to be induced by the
splitting  $(\D_V)_1 = \O_V \oplus T_V$ for any non-singular
variety $V$), so is a split almost polynomial algebra of
differential operators on $X_{i,a}$. 
\end{lemma}

\begin{definition}\rm 
For any complex scheme $T$, a {\bf family} 
$$\EE_T =(E_{i,T},\,s_{i,T},\,t_{i,T})$$
of pre-$\D$-modules parametrized
by $T$ is defined as follows. The $E_{i,T}$ are vector bundles
on $X_i\times T$ with structure of $\D_{i,T}$ modules, 
where $\D_{i,T}$ are the relative versions of the $\D_i$. 
The morphisms $s_{i,T}$ and $t_{i,T}$ are the relative versions
of the morphisms $s_i$ and $t_i$ in the definition of a
pre-$\D$-module. 
\end{definition}

Given any morphism $f:T'\to T$ of complex schemes and a family
$\EE_T$ of pre-$\D$-modules parametrized by $T$, the pullback
family $f^*\EE_T$ on $T'$ has again the obvious  
definition. This therefore defines a fibered category over the 
base category of complex schemes. When
we put the restriction that all morphisms in each fiber category
be isomorphisms, we get a fibered category $\P\!\D$ of groupoids over
$Schemes_{\C}$.

\begin{proposition}
The fibered category of groupoid $\P\!\D$ of pre-$\D$-modules on
$(X,Y)$ is an algebraic stack in the sense of Artin. 
\end{proposition}

\proof (Sketch) We refer the reader to the notes of Laumon [L]
for basic concepts and constructions involving algebraic stacks.
As fpqc descent and fpqc effective descent is
obviously satisfied by $\P\!\D$, it follows that $\P\!\D$ is a stack. It
remains to show that this stack is algebraic in the sense of
Artin. For this, first note that if $\Lambda$ is a split almost
polynomial algebra of differential operators, then $\O$-coherent
$\Lambda$-modules form an algebraic stack, for the forgetful
functor (1-morphism of stacks) from
$\Lambda$-modules to $\O$-modules is representable (as follows
from the alternative description of the structure of a
$\Lambda$-module given in lemma 2.13 of [S]), and
coherent $\O$-modules  form an algebraic stack in the sense of Artin.
Now from the lemma \ref{nnreplem} on coherence and representability of
the integrable direcrt image functor it can be seen that 
the forgetful functor (1-morphism of stacks) from pre-$\D$-modules to 
the product of the stacks of its underlying $\D_i$-modules is
representable (in fact, the details of this occur below in our
construction of a local universal family for pre-$\D$-modules).
Hence the result follows.

\subsection{Filtrations and 1-parameter deformations}
We first recall some standard deformation theory, for convinience
of reference. Recall that in this paper a vector bundle means a
locally free module (but not necessarily of constant rank), and
a subbundle of a vector bundle will 
mean a locally free submodule such that the quotient is also
locally free. 
Let $E$ be a vector bundle on a complex 
scheme $X$ together with an exhaustive increasing filtration $E_p$
by vector subbundles, indexed by $\Z$. (The phrase {\bf
exhaustive} means that $E_p=0$ for $p\ll 0$ and $E_p=E$
for $p\gg 0$.) Let $A^1 = \spec \C[\tau]$ be the affine line,
and let $U =\spec \C[\t, \t^{-1}]$ 
be the complement of the origin with inclusion
$j:U\hookrightarrow A^1$. Let $\pi_X:X\times A^1 \to X$ be the
projection. Consider the quasi-coherent sheaf 
$(1_X\times j)_*(\pi_X^*E|X\times U)$ on $X\times A^1$, which is
usually denoted by $E\otimes \C[\t,\t^{-1}]$. This has a
subsheaf $\ov{E}$ generated by all local sections of the
type $\t^p v_p$ where $v_p$ is a local section of
$\pi^*E_p$. It is common to write 
$$\ov{E} = \sum_{p\in\Z}E_p\tau^p\subset E\otimes
\C[\tau,\tau_{}^{-1}]$$
Then we have the following basic fact:

\rem\label{fact1} It can be seen that $\ov{E}$ is an
$\O_{X\times A^1}$-coherent submodule, which is in fact locally
free, and  $\ov{E}|X\times U$ is just $\pi_X^*E$ where
$\pi_X:X\times U\to X$ is the projection. On the
other hand, the specialization
of $\ov{E}$ at $\t =0$ is canonically isomorphic to the graded
object $E'=\oplus(E_p/E_{p-1})$ associated with $E$.  So, the 1-parameter
family $E_{\t}$ is a deformation of $E$ to its graded object $E'$.

Now let $F$ vector bundle on $X$ together with filtration $F_p$,
and let $f : E\to F$ be an $\O_X$-homomorphism. We 
have an induced $\O_{X\times A^1}$-homomorphism
$$\pi_X^*f\,:\,E\otimes \C[\t,\t^{-1} ]\to F\otimes \C[\t,\t^{-1}
]$$ 
Then we have the following basic fact :

\rem\label{fact2} 
The homomorphism $f:E\to F$ is filtered, that is, $f$ maps each $E_p$ into
$F_p$, if and only if the above homomorphism $\pi^*f$ carries
$\ov{E}$ into $\ov{F}$. In that case, the induced map at $\t=0$ is the
associated graded map $gr(f) : E'\to F'$.

As a consequence, we get the following :

\rem\label{fact3}
Let $E$ and $F$ be vector bundles with exhaustive filtrations
$E_p$ and $F_p$ indexed by $\Z$, and let $\ov{E}$ and
$\ov{F}$ be the corresponding deformations parametrized by
$A^1$. Let $f:E\to F$ be an $\O_X$-homomorphisms, and let 
$g:\ov{E}\to \ov{F}$ be an $\O_{X\times A^1}$-homomorphism.
Suppose that the restriction of $g$ to $X\times U$ (where
$U=A^1-\{ 0\}$) is equal to the pullback $\pi_X^*f$ of $f$ 
under $\pi_X:X\times U\to X$.
Then $f$ preserves the filtrations, and $g_0:E'\to F'$ is
the associated graded homomorphism $f':E'\to F'$  
(where $E'$ and $F'$ are the graded objects)

\begin{definition}\rm A {\bf sub pre-$\D$-module} $F$ of
of a pre-$\D$-module $E=(E_i,\,s_i,\,t_i)$ consists of the following
data : For each $i$ we are given an $\O_{X_i}$-coherent 
$\D_i$-submodule $F_i\subset E_i$ such that for each $i$,  

(i) the $\O_{X_i}$-modules $F_i$ and $E_i/F_i$ are
locally free (but not necessarily of constant ranks over $X_i$).
In other words, for each $(i,a)$, we are given a vector subbundle
$F_{i,a}\subset E_{i,a}$ which is a sub $\D_{i,a}$-module.

(ii) the maps $s_i$ and $t_i$ preserve $F_i$, that is 
$s_i$ maps $F_i|Y^*_i$ into $F_{i+1}|Y^*_i$ and $t_i$ maps
$F_{i+1}|Y^*_i$ into $F_i|Y^*_i$.

Note that consequently,  
$F=(F_i,\,s^F_i,\,t^F_i)$ is also a pre-$\D$ module, where $s^F$
and $t^F$ denote the restrictions of $s$ and $t$. Also, the
quotients $E_i/F_i$ naturally form a
pre-$\D$-module $E/F$ which we call as the corresponding {\bf
quotient pre-$\D$-module}.  
\end{definition}

\rem Given a pre-$\D$-module $\EE$ and a 
collection $\FF$ of subbundles $F_i\subset E_i$ which are
$\D_i$-submodules, 
the job of checking whether these subbundles are preserved by
the $s_i$ and $t_i$ is made easier by
the following: it is enough to check this in the fiber of a
point of each of the connected components of
$Y^*_i-p^{-1}(S_{i-1})$ where $p:Y^*_i\to S_i$ is the projection.
This is because the $s_i$ and $t_i$ are `integrable' in a suitable
sense, and if an integrable section $\sigma$ of a
vector bundle with an integrable connection has a value $\sigma(P)$ in
the fibre at $P$ of a subbundle preserved by the connection,
then it is a section of this subbundle.

\begin{definition}\rm  An {\bf exhaustive filtration on a
pre-$\D$-module} $\EE=(E_i,\,s_j,\,t_j)$ consists of an 
increasing sequence $\EE_p$ of sub pre-$\D$-modules of $\EE$ indexed by
$\Z$ such that $\EE_p=0$ for $p\ll 0$ and $\EE_p=\EE$
for $p\gg 0$.) A filtration is {\bf nontrivial} if $\EE_p$
is a nonzero proper sub pre-$\D$-module of $\EE$ for some $p$.  
\end{definition}

For a filtered pre-$\D$-module, each step
$\EE_p=((E_i)_p,\,(s_i)_p,\,(t_i)_p)$  of
the filtration, as well as the associated graded object
$\EE'=(E'_i,\,s'_i,\,t'_i)$ are pre-$\D$-modules.

Applying the remark \ref{fact2} to filtrations of
pre-$\D$-modules we get the following. 

\rem\label{deform}
Let $\EE$ be a pre-$\D$-module, together with an exhaustive
filtration $\EE_p$. Then
there exists a family $(\EE_\tau)_{\tau\in A^1}^{}$ of
pre-$\D$-modules parametrized by the affine
line $A^1=\spec\C[\tau]$, for which the specialization at
$\tau=0$ is the graded 
object $\EE'$ while the family over $\tau_0\neq0$ is the
constant family made from the original pre-$\D$-module 
$\EE$ defined as follows: put  
$\EE_i=\sum_{p\in\Z}(E_i)_p\tau^p\subset E_i\otimes
\C[\tau,\tau^{-1}]$.

\subsection{Quot scheme and group action on total family}
Let $X$ be a projective scheme over a base $S$ and $\V$ a coherent
$\O_X$-module. Let 
$$G =Aut_{\V}: {\rm Schemes}/S \to {\rm Groups}$$
be the contrafunctor which associates to any $T\to S$ the group
of all $\O_{X_T}$-linear automorphisms of the pullback $\V_T$ of
$\V$ under $X_T=X\times_ST\to X$. Note that then $G$ is in fact
an affine group scheme over $S$, but this will not be relevant
to us.

Let $Q=Quot_{\V/X/S}$ be the
relative quot scheme of quotients of $\V$ on fibers of $X\to S$.
A $T$-valued point $y:T\to Q$ is represented by a surjective 
$\O_{X_T}$-linear homomorphism $q:\V_T\to F$ where $F$ is a
coherent sheaf on $X_T$ which is flat over $T$. Two such 
surjections $q_1:\V_T\to F_1$ and $q_2:\V_T\to F_2$ represent
the same point $y\in Q(T)$ if and only if either of the
following two equivalent conditions is satisfied: (i) there
exists an isomorphism $f:F_1\to F_2$ such that $q_2=f\circ q_1$,
or (ii) the kernels of $q_1$ and $q_2$ are identical as a
subsheaf of $\V_T$. Therefore, a canonical way to represent the
point $y\in Q(T)$ is the quotient $\V_T\to \V_T/K_y$ where
$K_y=ker(q)$ depends only on $y$.  

A natural group action $Q\times G\to Q$ over $S$ is defined as
follows: 
in terms of valued points $y\in Q(T)$ represented by 
$q:\V_T\to F$, and $(g:\V_T\to \V_T) \in G(T)$, the point $yg\in Q(T)$ 
is represented by $q\circ g:\V_T\to F$. In other words, if $y$
is canonically represented by $\V_T\to \V_T/K_y$ then $yg$ is
canonically represented by $\V_T\to \V_T/g^{-1}(K_y)$. This means
$K_{yg}= g^{-1}(K_y)$.  

Let $q: \V_Q\to \U$ be the universal quotient family on $X_Q$.
The action $Q\times G\to Q$ over $S$ when pulled back under
$X\to S$ gives an action $Q_X \times G_X \to Q_X$ over $X$
(where $G_X=G\times _SX$ and $Q_X=Q\times_SX=X_Q$).  
The action $Q_X\times G_X\to Q_X$ has a natural lift to the sheaf $\U$
on $Q_X$ as follows. For $y\in Q(T)$ and $g\in G(T)$, the  pull
backs $\U_y$ and $\U_{yg}$ of the universal quotient sheaf under 
$y$ and $yg$ are canonically isomorphic to $\V_T/K_y$ and
$\V_T/g^{-1}(K_y)$ respectively, so $g:\V_T\to \V_T$ induces a canonical
isomorphism $\varphi_g^y: \U_{yg}\to \U_g$. Hence we get an
isomorphism $\varphi : g^*(\U)\to \U$ over $X_Q$. Since 
$\varphi_{gh}^y = \varphi_g^y\circ \varphi_h^{yg}$ for any 
$g,\,h\in G(T)$ and $t\in Q(T)$, we get 
$$\varphi_{gh} = \varphi_g\circ g^*(\varphi_h)$$
Thus, $\varphi$ is a `factor of automorphy', and so defines the
required lift.

\rem\label{liftformula}  
By definition, if $y\in Q(T)$ is represented by the surjection
$q:\V_T\to F$, then we get a canonical identification
of $F$ with $\U_y=\V_T/K_y$. For $g\in G(T)$ the point $yg\in Q(T)$ is
represented by $q\circ g:\V_T\to F$, so we get another canonical
identification 
of $F$ with $\U_{gy}=\V_T/g^{-1}(K_y)$. Under these
identifications, the isomorphism $\varphi_g^y:\U_{yg}\to U_y$
simply becomes the identity map $1_F :F\to F$. (This is so
simple that it can sometimes cause confusion.)
Hence if $\sigma$ is a local section of $\U_T$ 
represented by $(q,\, s)$, where $q:\V_T\to F$ and $s$ is a local
section of $F$, then the action of $g\in G(T)$ can be written as
$$(q,\,s)\cdot g = (q\circ g,\, s)$$

\rem\label{scalarmulti} 
The central subgroup scheme $\GG_m\subset G$ (where by 
definition $\lambda\in \GG_m(T)=\Gamma(T,\O_T^{\times})$ acts by scalar
multiplication on $\V_T$) acts trivially on $Q$ but its action on
the universal family $\U$ is again by scalar multiplication so
is non-trivial. In terms of the above notation have the equality
$$(q,\,s)\cdot \lambda = (q\circ \lambda, \,s) = (q, \,\lambda s)$$
Hence the induced action of $PG=G/\GG_m$ on $Q$ does not lift to $\U$,
which is the basic reason why a Poincar\'e bundle does not in general exist
in the kind of moduli problems we are interested in.

\rem\label{glnaction} 
In the applications below, $X$ will be in general a
projective scheme over $\C$ and the sheaf $\V$ will be of the
type $\O_X^{\oplus N}\otimes_{\O_X}\W = \C^N\otimes_{\C}\W$
where $\W$ is some coherent sheaf over $X$. Then $GL(N)$ is
naturally a subgroup scheme of $G=Aut_{\V}$, and we will only be
interested in the resulting action of $GL(N)$.

\subsection{Semistability and moduli for $\Lambda$-modules}

We now recall the moduli construction of Simpson for
$\Lambda$-modules (see section 2, 3 and 4 of [S] for details).
Let $X$ be projective, with ample line bundle $\O_X(1)$. 

Let $E$ be an $\O_X$-coherent $\Lambda$-module on $X$. Then 
Simpson defines $E$ to be a {\bf semistable $\Lambda$ module} if

(i) the $\O_X$-module $E$ is pure dimensional, and  

(ii) for each non-zero $\O_X$-coherent $\Lambda$-submodule
$F\subset E$, the inequality 
$$\dim\, H^0(X,\,F(n))/{\rm rank}\,(F)  \le 
{\rm dim}\,H^0(X,\,E(n))/{\rm rank}\,(E) $$
holds for $n$ sufficiently large (where ${\rm rank}\,(F)$ for any
coherent sheaf on $(X,\,\O_X(1))$ is by definition the leading
coefficient of the Hilbert polynomial of $F$). If the 
$\Lambda$-module $E$ is nonzero and moreover we can always have strict
inequality in the above for $0\ne F\ne E$, then $E$ 
is called {\bf stable}.

An {\bf S-filtration} of a semistable $\Lambda$-module $E$ is a
filtration $0=E_0\subset E_1\subset \ldots E_{\ell}=E$ by $\O_X$-coherent
$\Lambda$-submodules, such that each graded piece $E_i/E_{i-1}$
is a semistable $\Lambda$-module, with the same normalized
Hilbert polynomial as that of $E$ if non-zero
(where `normalized Hilbert polynomial' means Hilbert polynomial
divided by its leading coefficient). 
It can be seen that an S-filtration on a non-zero $E$ is maximal
(that is, can not be further refined) if and only if each graded piece
$E_i/E_{i-1}$ is stable. The associated graded
object $\oplus_{1\le i\le \ell}(E_i/E_{i-1})$ to a maximal 
S-filtration, after forgetting 
the gradation, is independent (upto isomorphism) of the choice of an S- 
filtration, and two non-zero semistable $\Lambda$-modules are called
{\bf S-equivalent} if they have S-filtrations with isomorphic
graded objects (after forgetting the gradation). The zero module
is defined to be S-equivalent to itself.

By a standard argument using Quot schemes (originally due to Narasimhan
and Ramanathan), it can be seen that semistability is a Zariski open
condition on the parameter scheme of any family of $\Lambda$-modules.

In order to construct a moduli for semi-stable
$\Lambda$-modules, Simpson first shows that if we fix the
Hilbert polynomial $P$, then all semi-stable
$\Lambda$-modules whose Hilbert polynomial is $P$ form a bounded
set, and then shows the following:

\begin{proposition}\label{simplocuni} 
(Simpson [S]) : Let $(X,\,\O_X(1))$ be a projective scheme, $P$
a fixed Hilbert polynomial, and $\Lambda$ a sheaf of
differential operators on $X$.
Then there exists a quasi-projective scheme $C$
together with an action of $PGL(N)$ (for some large $N$), and a
family $E_C$ of semistable $\Lambda$-modules on $X$ with Hilbert polynomial
$P$ para\-met\-rized by $C$ such that 

(1) the family $E_C$ is a local universal family for semistable
$\Lambda$-modules with Hilbert polynomial $P$,

(2) two morphisms $f_1,f_2:T\stackrel{\to}{\to}C$ give
isomorphic families $f_1^*(E_C)$ and $f_2^*(E_C)$ if and only if
there exists a Zariski open cover $T'\to T$ (that is, $T'$ is a
disjoint union of finitely many open subsets of $T$ whose union
is $T$) and a $T'$-valued point $g:T'\to PGL(N)$ of $PGL(N)$
which carries $f_1|T'$ to $f_2|T'$. 

(3) a good quotient (`good' in the technical sense of GIT) $C\to
C/\!/PGL(N)$ exists, and is (as a consequence of (1) and (2)) the
coarse moduli scheme for $S$-equivalence classes of semistable
$\Lambda$-modules with Hilbert polynomial $P$.

(4) {\bf Limit points of orbits }: Let $\lambda :GL(1)\to GL(N)$ be a
1-parameter subgroup, and let $q\in C$ have limit $q_0\in C$
under $\lambda$, that is,
$$q_0 =\lim_{\t\to 0}\, q\cdot \lambda(\t)$$
Let $E_{\t}$ be the pullback of $E_C$ to $X\times A^1$ under the
resulting morphism $\lambda :A^1\to C$ on the affine line $A^1$.
Then there exists an exhaustive filtration $F_p$ of $E$ (where
$E$ is the module associated to $\t =1$) by $\Lambda$-submodules
such that each $F_p/F_{p-1}$ is semistable, and when non-zero it has
the same reduced Hilbert polynomial as $E$, 
and the family $E_{\t}$ is isomorphic to the family
$\ov{E}=\sum_{p\in \Z} F_p \t^p \,\subset \, 
E \otimes \C[\t,\, \t^{-1}]$ occuring in remark \ref{fact1} above. 
In particular, the limit $E_0$ is the graded object corresponding to
$F_p$, and so the closed points of the quotient are in bijection with
the set of $S$-equivalence classes. 
\end{proposition}

\rem 
The following feature of the family $E_C$ and the action of
$GL(N)$ will be very important: $C$ is a certain locally closed
subscheme of the Quot scheme $Quot_{\C^N\otimes \W/X/\C}$ of quotients   
$$q :\C^N\otimes \W \to E$$
where $W$ is some fixed coherent sheaf on $X$, 
such that $C$ is invariant under $GL(N)$, and $E_C$ is the
restriction to $C$ of the universal family $\U$ on the Quot scheme,
The action of $GL(N)$ on the Quot scheme and its lift to $\U$ is as
explained in remark \ref{glnaction} above.

\vfill
\pagebreak
\subsection{Strong local freeness for semistable $\Lambda$-modules}

\begin{definition}\label{defstrong}\rm
We will say that a semistable $\Lambda$-module $E$ on $X$ is {\bf
strongly locally free} if for every S-filtration 
$0=E_0\subset E_1\subset \ldots E_{\ell}=E$, the associated
graded object $\oplus_{1\le i\le \ell}(E_i/E_{i-1})$ is a locally
free $\O_X$-module. (In particular, the zero module is strongly
locally free.) 
\end{definition}

\rems (1)  A nonzero semistable $\Lambda$-module $E$ is strongly
locally free if and only if there exists some maximal S-filtration such
that the graded pieces are are locally free $\O_X$-modules. 

(2) A strongly locally free semistable $\Lambda$-module is
necessarily locally free. However, a locally free semistable
$\Lambda$-module is not necessarily strongly locally free.

\begin{proposition}\label{strong}  
(1) Let $E_C$ be the family of semistable $\Lambda$-modules
parametrized by $C$. Then the condition that the associated
graded object to any S-filtration $F_k$ of $E_q$ should be
locally free defines a $GL(N)$-invariant open subset  
$C^o$ of $C$ which is closed under limits of orbits, and hence has a
good quotient under $GL(N)$, which is an open subscheme of the
moduli of semistable $\Lambda$-modules.

(2) For any family $E_T$ of semistable
$\Lambda$-modules parametrized by a scheme $T$, the condition
that $E_t$ is strongly locally free is a Zariski open condition
on $T$.
\end{proposition} 

\proof Let $U\subset C$ be the open subscheme defined by the
condition that $E_q$ is locally free for $q\in U$ (this is
indeed open as a consequence of the flatness of $E_C$ over $C$).
Note that $U$ is $G=GL(N)$-invariant. Let $\pi :C\to C/\!/G$ be
the good quotient. Then as $F=C-U$ is a $G$-invariant closed
subset, $\pi(F)$ is closed in $C/\!/G$ by properties of good
quotient. Let $C^o = \pi^{-1} (C/\!/G -\pi(F))$. Then $C^o$ is the
desired open subscheme of $C$, which proves the statement (1). 

It follows from the relation between closure of orbits
and S-filtrations that points of $C/\!/G$ correspond
to S-equivalence classes. Therefore the statement (2) now
follows from (1) by the universal property of the moduli $C/\!/G$.

\subsection{Semistable pre-$\D$-modules --- definition}

We now come back to $(X,Y)$ as before, and our earlier notation.
We choose an ample line bundle on each $X_{i,a}$, and 
fix the resulting Hilbert polynomials
$p_{i,a}(n)$ of the sheaf $E_{i,a}$ on $X_{i,a}$. 
By lemma \ref{lambda}  the theory of semistability and moduli for
$\Lambda$-modules can be applied to $\D_{i,a}$-modules.  

\begin{definition}\rm We say that a pre-$\D$-module is {\bf
semistable} if the following two conditions are satisfied.

(i) Each of the $\D_{i,a}$-modules 
$E_{i,a}$ is semistable in the sense of Simpson. 

(ii) Each $E_{i,a}$ is strongly locally free, that is,  
the associated graded object to $E_{i,a}$ under any S-filtration
is again a locally free $\O_{X_{i,a}}$-module (see definition
\ref{defstrong} and its subsequent remarks).  
\end{definition}

\rem The notion of stability is treated later.

\begin{proposition}
A vector bundle $E_{i,a}$ on $X_{i,a}$ with the structure of a
semistable $\D_{i,a}$-module is strongly locally free if it has 
good residual eigenvalues as defined in definition
\ref{goodresieigen}.  

As a consequence, a pre-$\D$-module $E$ such that each of the
$\D_{i,a}$-modules $E_{i,a}$ is semistable and has good residual
eigenvalues is a semistable pre-$\D$-module. 
\end{proposition}

\proof Let $Z$ be the polydisk in $\C^n$ defined by $|z_i|< 1$, and
$W\subset Z$ the divisor with normal crossings
defined by $\prod_{k\le r}z_k =0$ (if $r=0$ then $W$ is empty). 
Let there be given
a set theoretic section (fundamental domain) for the exponential
map $\C\to \C\,^* \,:\,z\mapsto \exp(2\pi iz)$. Then the Deligne
construction, which associates to a local system $\L$ on $Z-W$ an
integrable logarithmic connection $E(\L)$ on $(Z,W)$ with residual
eigenvalues in the given fundamental domain, has the following
property: if $\K\subset\L$ are two local systems, then $E(\K)$ is a
subbundle of $E(\L)$. Now suppose $F$ is an $\O_Z$-coherent
$\D_Z[\log W]$-submodule of $E(\L)$. Let $\K$ be the local
system on $Z-W$ defined by $F$, and let $E(\K)$ be its
associated logarithmic connection given by Deligne's
construction. Then $F|(Z-W)=E(\K)|(Z-W)$. As $E(\K)$ is a vector
subbundle of the vector bundle $E(\L)$, it follows that 
$E(\K)$ is just the $\O_Z$-saturation of $F$. Hence if
$E$ has good residual eigenvalues and if $F\subset E$ is an
$\O_Z$-coherent and saturated $\D_Z[\log W]$-submodule, then $F$
is a vector subbundle. 

Now we apply this to each $E_{i,a}$ as follows. Let 
$E=E_{i,a}$ be $\D_{i,a}$-semistable, and let $F\subset E$
be a step in an S-filtration of $E$. 
Hence we must have 
$$ {p_F\over rank(F)} =  {p_E\over rank(E)}$$ 
where $p_F$ denotes the Hilbert polynomial of a sheaf $F$.
Note that the
$\O_{X_{i,a}}$-saturation $F'$ of $F$ in $E$ is again an
$\O$-coherent $\D_{i,a}$-submodule, which restricts to $F$ on a
dense open subset of $X_{i,a}$, in particular, $F'$ has the same rank
as $F$. Hence with respect to any ample line bundle on
$X_{i,a}$, the normalized Hilbert
polynomials of $F$ and $F'$ satisfy the relation
$$ {p_F\over rank(F)} \le  {p_{F'}\over rank(F')}$$ 
with equality only if $F=F'$. 
As  $E=E_{i,a}$ is semistable, we have 
$${p_{F'}\over rank(F')} \le {p_{E}\over rank(E)}$$ 
and hence $F=F'$. Hence we can assume that any step $F$ in an
S-filtration of $E=E_{i,a}$ is $\O$-saturated.

Hence the proposition would follow if we show that if for each
$i\ge m+1$, $E_{i,a}$
has good residual eigenvalues under $\theta_{i-1}$ (see
definition \ref{goodresieigen}), 
then any $\O$-coherent and saturated $\D_{i,a}$-submodule $F$ is
a vector subbundle. This is a purely local question on $X$ in
the euclidean topology, so we may assume that $X$ is a polydisk
with local coordinates $z_i$ and $Y$ is defined by a monomial in
the $z_i$. Then $X_{i,a}$ is just some coordinate $i$-plane $Z$
contained in $Y$, and $W=X_{i,a}\cap S_{i-1}$ is a normal
crossing divisor in the polydisk $Z$, defined by a monomial in
the coordinates. Using the local coordinates $z_k$, we get a
structure of $\D_Z [\log W]$-module on the vector bundle
$E=E_{i,a}$ on $Z$ with good residual eigenvalues, and $F$
becomes an $\O_Z$-coherent and saturated $\D_Z[\log
W]$-submodule. (Note that these $\D_Z[\log W]$-structures very
much depend on the choice of local coordinates $z_k$). As $E$ has good
residual eigenvalues by assumption, it follows from the first
part of our argument (involving Deligne constructions over
polydisks) that $F$ is a vector subbundle of $E_{i,a}$. This
proves the proposition.  

\rem\label{goodresi2} If $E_{i,a}$ has good residual eigenvalues under
$\theta_{i-1}$, and if $F$ is a vector subbundle which is a
$\D_i$-submodule, then it can be seen that the associated graded
$\D_i$-module $F\oplus (E_{i,a}/F)$ also has good residual
eigenvalues under $\theta_{i-1}$.

The following example shows that if $E$ is a vector bundle with
an integrable logarithmic connection, and $F\subset E$ is an
$\O$-coherent sub connection such that $E/F$ is torsion free,
then it can still happen that $E/F$ is not locally free if $E$
has bad residual eigenvalues.

\begin{example}\label{Esnault}\rm (Due to H\'el\`ene Esnault)
Let $X$ be a 
polydisk in $\C^2$, with divisor $Y$ defined by $xy=0$. Let
$E=\O_X^{\oplus 2}$ with basis $v_1=(1,0)$ and $v_2=(0,1)$. 
On $E$ we define an integrable logarithmic connection 
$\nabla :E\to \Omega^1_X[\log Y]\otimes E$ by putting
$\nabla(v_1) = (dx/x)\otimes v_1 $ and
$\nabla(v_2)=(dy/y)\otimes v_2$ (this has curvature zero, as
$(-\log x)v_1$ and $(-\log y)v_2$ form a flat basis of $E|X-Y$). 
Note that both $0$ and $1$ are eigenvalues of the residue of
$(E,\nabla)$, along any branch of $Y$.
Let ${\bf m}\subset \O_X$ be the ideal sheaf generated by $x$ and $y$.
This is a sub $\D_X[\log Y]$-module of $\O_X$, which is torsion
free but not locally free. Now define a surjective homomorphism
$\varphi : E\to {\bf m}$ sending $v_1\mapsto x$ and $v_2\mapsto y$,
which can be checked to be $\D_X[\log Y]$-linear, and put
$F=ker(\varphi)$. Then $F$ is an $\O_X$-saturated subconnection
(which is in fact locally free), but $E/F={\bf m}$ is not
locally free.  
\end{example}

\subsection{Semistable pre-$\D$-modules --- local universal family}

We now construct a local universal family for semistable
pre-$\D$-modules with given Hilbert polynomials. Let $C_{i,a}$
be the scheme for $(\D_{i,a},\,p_{i,a})$ given by the
proposition \ref{simplocuni}, with the action of
$PGL(p_{i,a}(N_{i,a}))$ as in the proposition \ref{simplocuni}.
Let $C^o_{i,a}\subset C_{i,a}$ be the open subset where
$E_{i,a}$ is strongly locally free.  By proposition
\ref{strong}, $C^o_{i,a}$ is an open subset of $C_{i,a}$ which
is $PGL(p_{i,a}(N_{i,a}))$-invariant and admits a good quotient
for the action of $PGL(p_{i,a}(N_{i,a}))$, which is an open
subscheme of $C_{i,a}/\!/PGL(p_{i,a}(N_{i,a}))$.

Let $C^o_i= \prod_a C^o_{i,a}$ and let 
$C=\prod_i C^o_i=\prod_{i,a} C^o_{i,a}$. Let $E_{i,a}$ again
denote the pullback of $E_{i,a}$ to $X_{i,a}\times C$ under the
projection $C\to C^o_{i,a}$. Let $E_i|Y^*_i$ and $E_{i+1}|Y^*_i$
be regarded as families of $\D^*_i$-modules parametrized by $C$. 
These are again flat over $C$ as the $E_i$ are locally free on
$X_i$. Hence by applying lemma 
\ref{nnreplem} to the pair $E_{i+1}|Y^*_i$ and
$E_i|Y^*_i$ of $\D^*_i$-modules parametrized by $C$, we
get linear schemes $A_i$ and $B_i$ over $C$ which parametrize
$\D^*_i$-linear maps  $t_i$ and $s_i$ in either direction
between the specializations of these two families. Let
$H_i\subset A_i\times_CB_i$ be the closed subscheme defined by
the conditions on $t_i$ and $s_i$ imposed by the definition of a
pre-$\D$-module (it can be seen that the conditions indeed
define a closed subscheme $H_i$). 
Finally, let $H$ be the fibered product over $C$ of all the
$H_i$. By its construction, $H$ parametrizes a natural family of
pre-$\D$-modules on $(X,Y)$.

Let $G_i =\prod_a G_{i,a}$ and let 
$G = \prod_i G_i=\prod_{i,a} G_{i,a}$. Note that this is a
reductive group. We define an action of $G$ on $H$ as follows.
Any point $\un{q}_{i,a}$ of $C^o_{i,a}$ is represented by a quotient 
$q_{i,a} :\O_{X_{i,a}}(-N_{i,a})^{p_{i,a}(N_{i,a})} \to  E_{i,a}$ 
(which satisfies some additional properties) and a point of $H$ over a
point $(\un{q}_{i,a})\in C=\prod_{i,a}C^o_{i,a}$ is given by the
additional data $s_j:(E_j|Y^*_j)\to (E_{j+1}|Y^*_j)$ and 
$t_j:(E_{j+1}|Y^*_j)\to (E_j|Y^*_j)$, and
so a point of $H$ is represented by the data $(q_i,s_j,t_j)$. 

\rem\label{eq} Note that two such tuples $(q_i,s_j,t_j)$ and
$(q'_i,s'_j,t'_j)$ 
represent the same point of $H$ if there exists an isomorphism 
$\phi: E \to E'$ of pre-$\D$-modules $E=(E_i,s_it_i)$ and 
$E' =(E'_i,s'_i,t'_i)$ such that $q'_i = \phi_i\circ q_i$ for
each $i$. 

\begin{definition}\label{defaction}
\rm (Right action of the group $G$ on the scheme $H$.)
In terms of valued points, we define this as follows.
For any point $h$ of $H$ represented by $(q_i,s_j,t_j)$, and
an element $g=(g_i)\in G = \prod_i G_i$, put  
$$ (q_i,\,s_j,\,t_j)\cdot g = (q_i\circ g_i,\,s_j,\,t_j )$$
Note that this is well defined with respect to the equivalence
given by remark \ref{eq}, and indeed defines an action of $G$ on
$H$ lifting its action on $C$, as follows from remark
\ref{liftformula}. 
\end{definition}

It is clear from the definitions of $H$ and this action that two
points of $H$ parametrise isomorphic pre-$\D$-modules if and only
if they lie in the same $G$ orbit. 

The morphism $H\to C$ is an affine morphism
which is $G$-equivariant, where $G$ acts on $C$ via $\prod_{i,a}
G_{i,a}$. As seen before, the action of each $G_{i,a}$ on $C^o_{i,a}$
admits a good quotient in the sense of geometric 
invariant theory, and hence the action of $G$ on $C$ admits a
good quotient $C/\!/G$. A well known lemma of Ramanathan
(see Proposition 3.12 in [Ne]) asserts that if $G$ is a
reductive group acting on two schemes $U$ and $V$ such that $V$
admits a good quotient $V/\!/G$, and if there exists an affine,
$G$-equivariant morphism $U\to V$, then there exists a good
quotient $U/\!/G$. Applying this to the $G$-equivariant affine
morphism $H\to C$, a good quotient $H/\!/G$ exists, 
which by construction and universal properties of good quotients
is the coarse moduli scheme of semistable pre-$\D$-modules with
given Hilbert polynomials. By construction this is a separated
scheme of finite type over $\C$, and is, in fact, quasiprojective.  

Note that under a good quotient in the sense of geometric
invariant theory, two different orbits can in some cases get mapped
to the same point (get identified in the quotient). 
In the rest of this section, we determine what are the closed
points of the quotient $H/\!/G$. 

\subsection{Stability and points of the moduli}

Let $T_{i,a}\subset G_{i,a}$ be its center. Let
$T_i=\prod_aT_{i,a}$, and $T= \prod_iT_i=(\prod_{i,a}T_{i,a})$. 
Note that $T \subset G$ as a closed normal subgroup 
(which is a torus), which acts trivially on $C$. 
By definition of $T$ we have a canonical identification 
$$T = \prod_i \Gamma(X_i,\O^{\times}_{X_i})$$
By definition \ref{defaction} and remark \ref{scalarmulti}, the
action of $\lambda =(\lambda_i)\in T$ on 
$h=(q_i,s_i,t_i)\in H$ is given by 
$$h\cdot g =( \lambda_iq_i,\,s_i,\,t_i) = 
(q_i,\, {\lambda_i\over\lambda_{i+1}}s_i,\, 
{\lambda_{i+1}\over\lambda_i}t_i)$$
In [N-S],the construction of the quotient
$H/\!/G$ was made in a complicated way in two steps: 
by Ramanathan's lemma, we can first have the
quotient $R=H/\!/T$, and then as the second step we
have the quotient $H/\!/G = R/\!/(G/T)$. 
However, we now do it in a much simplified way, which gives a
simplification also of [N-S].

\begin{definition}\rm
A sub pre-$\D$-module $\FF$ of a semistable pre-$\D$-module $\EE$
will be called  an {\bf S-submodule} if each non-zero $F_{i,a}$
has the same normalized Hilbert polynomial $p_{i,a}$ as that of
$E_{i,a}$.  
A filtration $\EE_p$ on a semistable pre-$\D$-module $\EE$ is an 
{\bf S-filtration of the pre-$\D$-module} if each
$\EE_p$ is an S-submodule, equivalently the given
filtration on each $E_{i,a}$ is an S-filtration. 
\end{definition}

\rem Given a semistable pre-$\D$-module $\EE$, an S-submodule, the
corresponding quotient module, and 
the graded pre-$\D$-module $\EE'$
associated with an S-filtration are again semistable pre-$\D$-modules. 
Moreover, $\EE'$ has the same Hilbert polynomials $p_{i,a}$ as $\EE$.
We will say that $\EE$ is (primitively) S-equivalent to $\EE'$. 

\begin{definition}\rm
The equivalence relation on the set of isomorphism classes
of all semistable pre-$\D$-modules generated by the above relation,
under which the graded module $\EE '$ associated to an
S-filtration of $\EE$ is taken to be equivalent to $\EE$, will
be called {\bf S-equivalence for pre-$\D$-modules}.
\end{definition}

\begin{definition}\label{defstable}\rm
We say that a semistable pre-$\D$-module is {\bf stable} if it
nonzero and does not admit any nonzero proper S-submodule. 
\end{definition}

\begin{proposition}\label{1parafamily}
Let $\EE_H$ denote the tautological family of pre-$\D$-modules
para\-met\-rized by $H$. Let $\lambda :GL(1)\to G$ be a 1-parameter
subgroup of $G=\prod G_i$, and let $h=(q_i, s_j, t_j)\in H$ be a
point such that the limit $\lim_{\t \to 0}\,h\lambda(\t)$ exists
in $H$. Let $\ov{\lambda}:A^1\to H$ be the resulting morphism.
Then there exists an S-filtration $(\EE_h)p$ of the
pre-$\D$-module $(\EE_h)_p$ such that the pullback of $\EE_H$ to $A^1$
under $\ov{\lambda}:A^1\to H$ is isomorphic to the  
family constructed in remark \ref{deform}
\end{proposition}

\proof By the definition of the action of $G$ on $H$, the family
$\EE_{\tau}$ satisfies the following properties:

(i) The families $(E_{i,a})_{\t}$ are of the necessary type by 
proposition \ref{simplocuni}.(4).

(ii) Outside $\tau=0$, the homomorphisms $(s_j)_{\tau}$ and
$(t_j)_{\t}$ are pull backs of $s_j$ and $t_j$. 

Therefore now the proposition follows from remark \ref{fact3}.

The following lemma, whose proof is obvious, is necessary to
show that stability in an open condition on the parameter scheme
$T$ of a family $\EE_T$ of pre-$\D$-modules.

\begin{lemma}
Let $\EE_T$ be a family of pre-$\D$-modules parametrized by $T$.
Let there be given a family $\FF_T$ of sub vectorbundles
$F_{i,T}\subset E_{i,T}$ which are $\D_{i,T}$-submodules. Then
there exists a closed subscheme $T_o\subset T$ with the
following universal property. Given any base change $T'\to T$,
the pullback $\FF_{T'}$ is a sub pre-$\D$-module of $\EE_{T'}$
if and only if $T'\to T$ factors through $T_o$.
\end{lemma}

Using the above lemma, the `quot scheme argument for openness of
stability' can now be applied to a
family of pre-$\D$-modules, to give

\begin{proposition}\label{open}
Stability is a Zariski open condition on the
parameter scheme $T$ of any family $\EE_T$ of semistable
pre-$\D$-modules. 
\end{proposition}

\proof Let $\pi:P\to T$ be the projective scheme, which is closed
subscheme of a fibered product over $T$ of relative quot schemes
of the $E_{i,T}$, 
which parametrizes families of sub vector bundles
$F_i$ which are $\D_i$-submodules with the same reduced
Hilbert polynomials as those of $\EE$. (Actually, the Quot
scheme parametrizes coherent quotients flat over the base, and
the condition of $\D_i$-linearity gives a closed subscheme. Now
by the assumption of strong local freeness on the $E_{i,a}$, it
follows that the quotients are locally free). Now by the above
lemma, $P$ has a closed subscheme $P_o$ where $\FF_{P_o}$ is a
family of sub pre-$\D$-modules with the same normalized Hilbert
polynomials, and every such sub pre-$\D$-module of a $\EE_t$ for
$t\in T$ occures among these. Hence $T-\pi(P_o)$ is the desired
opoen subset of $T$. 

\rem As semistability is itself a Zariski open condition on any
family of pre-$D$-modules, it now follows that 
stability is a Zariski open condition on the
parameter scheme $T$ of any family $\EE_T$ of pre-$\D$-modules. 

Now all the ingredients are in place for the following main
theorem, generalizing theorem 4.19 in [N-S].

\begin{theorem}\label{maintheorem}  
Let $X$ be a non-singular variety with a normal crossing divisor
$Y$. Let a numerical polynomial $p_{i,a}$ and an ample line
bundle on $X_{i,a}$ be chosen for each $X_{i,a}$. Then we have
the following.  

(1) There exists a coarse moduli scheme $\M$ for semistable
pre-$\D$-modules $\EE$ on $(X,Y)$ where $E_{i,a}$ has Hilbert
polynomial $p_{i,a}$. The scheme $\M$ is quasiprojective, in
particular, separated and of finite type over $\C$.   

(2) The points of $\M$ are S-equivalence classes of semistable
pre-$\D$-modules.  

(3) The S-equivalence class of a stable pre-$\D$-module equals
its isomorphism class.  

(4) $\M$ has an open subscheme $\M ^s$ whose points are the
isomorphism classes of all stable pre-$\D$-modules. This is a
coarse moduli for (isomorphism classes of) stable pre-$\D$-modules. 
\end{theorem}

\proof The statement
(1) is by the construction of $\M=H/\!/G$ and properties of a
good quotient. 

The statement (2) follows from remark \ref{deform} and
propositon \ref{1parafamily}.  

Let $x\in H$ and let $x_0 \in H$ be a limit point of the
orbit $Gx$. Then by properties of GIT quotients, there exists a 
$1$-parameter subgroup 
$\lambda:GL(1)\to G$ such that $x_0=\lim_{\t\to
0}\,x\cdot\lambda(\t)$. Any such limit is of the type given by
proposition \ref{1parafamily}, made from an S-filtration of the
corresponding pre-$\D$-module $\EE_x$. 
If $x\in H$ be stable (means corresponds to a stable
pre-$\D$-module), then it has no non-trivial S-filtration, so
the orbit of $x$ is closed. As stability is an open condition 
on $H$ by proposition \ref{open}, if the orbit of a point $y$ in
$H$ has a limit $x$ which is stable, then the point $y$ must
itself be stable. So by the above, the orbit of $y$ must be
closed, so $x\in Gy$. Hence a stable point $x$ is not the limit
point of any other orbit.  Hence (3) follows. 

Finally, the statement (4) follows from (1), (2), (3), and
proposition \ref{open}. This completes the proof of the theorem.

{\footnotesize
Note: the statement (3) in theorem 4.19 of
[N-S] has a mistake - the `if and only if' should be changed to
`if', removing the `only if' part.) 

}

\bigskip

\section{Perverse Sheaves on $(X,Y)$}
In this section, we give a finite description (in terms of a
finite quiver of finite dimensional vector spaces and linear
maaps) of perverse
sheaves on $(X,Y)$, that is, perverse sheaves on $X$ that are
cohomologically constructible with respect to the stratification
$X=\cup_i\, (S_i-S_{i-1})$, which closely parallels our
definition of pre-$\D$-modules. This enables us to describe
the perverse sheaf associated to the $\D$-module associated to a
pre-$\D$-module directly in terms of the pre-$\D$-module. In
turn, this allows us to deduce properties of the analytic
morphism from the moduli of pre-$\D$-modules to the moduli of
such quivers (which is the moduli of perverse sheaves with given
kind of singularities).

More general finite descriptions of
perverse sheaves in terms of quivers exist in literature (see
for example MacPherson and Vilonen [M-V]), where the requirement
of normal crossing singularities is not needed, and the
resulting moduli space has been constructed by Brylinski,
MacPherson, and Vilonen [B-M-F]). We cannot use their moduli
directly, as it is too general for the specific purpose of
describing the Riemann-Hilbert morphism by a useable formula
(where, for example, we can see the differential of the map).

\subsection{The specialization functor}
In this section, a vector bundle will usually mean a geometric
vector bundle in the analytic category. That is, if $E$ is a locally
free sheaf (but not necessarily of constant rank) on a reduced
scheme $X$ of finite type over $\C$ , then when we refer to the
{\bf vector bundle $E$}, what we mean is the analytic space
(with euclidean topology) associated to the scheme 
$\spec_X Sym^{\cdot}(E^*)$.  

The replacement by Verdier (see [V1] and [V2]) of the not so
canonical operation of {\sl restriction to a tubular
neighbourhood} by {\sl specialization to normal cone (or
bundle)} is used below in a somewhat more general set up as
follows. 

Let $M$ be a complex manifold, $T \subset M$ be a divisor
with normal crossings, and for some integer $k$ 
let $T_k$ be a union of components of the
$k$-dimensional singularity 
stratum of $T$. Let $C \to T_k$ be a union of components of the
normalization of $T_k$, 
and $f: C\to M$ the composite map $C\to T_k\to M$. Let $N_f$ be
the normal bundle to $f:C\to M$, and let $U\subset N_f$ be the
open subset which is the complement of the normal crossing
divisor $F_f$ in $N_f$ defined by vectors tangent to branches of $T$.
Note that in particular $F_f$ contains the zero section of
$N_f$. Then we have a functor from local systems on $M-T$
to local systems on $U$ defined as follow.

For each $x\in C$, there exists an open neighbourhood $V_x$ of
$f(x)$ in $M$ such that the restricted map $f_x:C_x \to M$ where
$C_x=f^{-1}(V_x)$ and $f_x=f|C_x$ is 
is a closed imbedding of the manifold $C_x$ into
$V_x$. Then note that the normal bundle $N_{f_x}$ of $f_x:C_x\to
V_x$ is the restriction of $N_f$ to $C_x\subset C$. Let
$F_{f_x}=F_f\cap N_{f_x}$, and $U_x=N_{f_x}-F_{f_x}$. 
The usual functor of specialization (see [V1] and [V2]) now
associates a local system on $U_x$ to a local system on $V_x-T$.
These glue together to define our desired functor.
Given a local system $\E$ on $M_T$, we denote its specialization
by $\E|\!|U$, which is a local system on $U$. 

More generally, the above method gives a definition of a
specialization functor \\
between the derived categories of
cohomologically bounded constructible complexes of sheaves of
complex vector spaces on $M$ and $N_f$. If the complex
$\F^{\cdot}$ is cohomologically constructible with respect to
the singularity stratification of $(M, T)$, then its
specialization $\F^{\cdot}|\!|N_f$ to $N_f$ is 
cohomologically constructible with respect to
the singularity stratification of $(N_f,\,F_f)$. This functor
carries perverse sheaves to perverse sheaves.  

\rem\label{composite} For a topological manifold $M$ which is
possibly disconnected, 
choose a base point in each component $M_a$, and let $\Gamma^M_a$
denote the indexed set of the fundamental groups of the
components of $M$ with respect to the chosen base points,
indexed by $a\in \pi_0(M)$. We have an equivalence categories
between local systems on $T$ and an indexed collection of group
representations $\rho_a: \Gamma^T_a \to GL(n_a)$. 
In the above situation, we have the groups $\Gamma^U_a$ and
$\Gamma^{M-T}_b$. Then there exists a map $\gamma: \pi_0(U)\to
\pi_0(M-T)$ and a group homomorphism
$$\psi_a:\Gamma^U_a\to \Gamma^{M-T}_{\gamma(a)}$$
for each $a\in \pi_0(U)$, such that the above functor of
specialization from 
local systems on $M-T$ to local systems on $U$ is given by
associating to a collection of representations
$\rho_b :\Gamma^{M-T}_b \to GL(n_b)$ where $b$ varies over
$\pi_0(M-T)$ the collection composite 
representations $\rho_{\gamma(a)}\circ \psi_a$. In this sense,
specialization is like pullback.

\subsection{Finite representation}
We now return to $(X,Y)$ and use our standard notation (see
section 2). We will apply the above specialization functor to
the following cases.  

Case(1): $M$ is $X$, $T$ is $Y$, and $k=d-1$, so $C=X_{d-1}=Y^*$ is
the normalization of $Y$. In this case, 
starting from a local system $\E_d$ on $X-Y$ we get a local system on
$U_{d-1}$, which we denote by $\E_d|\!|U_{d-1}$.

Case(2): $M$ is $N_{i+1}$ ($=$ the normal bundle to $f_i:X_i\to X$), 
$T$ is the divisor $F_{i+1}$ in $N_{i+1}$ defined by vectors tangent to
branches of $Y$, $k=i$, $T$ the inverse image of $S_i$ under the
map $p_{i+1}:X_{i+1}\to S_{i+1}$, and $C=Y^*_i$.  
In this case, starting from a local system $\E_{i+1}$ on $U_{i+1}$,
we get a local system on $R_i$ which we denote by $\E_{i+1}|\!|R_i$.  

Case(3): Apply this with $M=N_{i+2}$, $T=F_{i+2}$, 
$C=Z_i$ and $C\to M$ the composite 
$$Z_i\stackrel{p_{\{i+2 \},\{ i,i+2\} }}{\to} X_{i+2}
\hookrightarrow N_{i+2}$$ 
In this case, starting from a local system $\E_{i+2}$ on $U_{i+2}$,
we get a local system on $W_i$ which we denote by $\E_{i+2}|\!|W_i$.  

Case(4): Apply this with $M=N_{i+1}$, $T=F_{i+1}$, 
$C=Z^*_i$ and $C\to M$ the composite 
$$Z^*_i\stackrel{p_{\{i+1 \},\{ i,i+1, i+2\} }}{\to} X_{i+1}
\hookrightarrow N_{i+1}$$ 
In this case, starting from a local system $\E_{i+1}$ on $U_{i+1}$,
we get a local system on $W^*_i$ which we denote by $\E_{i+1}|\!|W^*_i$.

Case(5): Apply this with $M=N_{i+1,i+2 }$, $T=F_{i+1,i+2}$, 
$C = Z^*_i$, and $C\to M$ the composite 
$$Z^*_i\stackrel{p_{\{i+1, i+2 \},\{ i,i+1, i+2\} }}{\to} Y^*_{i+1}
\hookrightarrow N_{i+1, i+2}$$ 
In this case, starting from a local system $\F_{i+1}$ on $R_{i+1}$,
we get a local system on $W^*_i$ which we denote by $\F_{i+1}|\!|W^*_i$.  

On the other hand, note that the derivative of the covering projection
$p_{\{i \},\{ i,i+1\} } : Y^*_i\to X_i$ is a map 
$dp : N_{i,i+1} \to N_i$ under which $F_{i,i+1}\subset
N_{i,i+1}$ is the inverse image of $F_i\subset N_i$. 
Hence $dp$ induces a map $R_i\to U_i$. If $\E_i$ is a local
system on $U_i$, then we denote its pullback under this map by
$\E_i|R_i$, which is a local system on $R_i$. 

Similarly, for the covering projection 
$p_{\{i \},\{ i,i+2\} } :Z_i\to X_i$ the derivative induces a
map $W_i\to U_i$. If $\E_i$ is a local system on $U_i$ then we
denote its pullback under this map by $\E_i|W_i$ which is a
local systems on $W_i$.  

Note that the derivative of the $2$-sheeted covering projection
$X_{\{ i,i+1,i+2 \} }\to X_{\{ i,i+2\} }$ induces a $2$ sheeted
covering projection $\pi : W^*_i\to W_i$. We will denote the
pullbacks of $\E_i|W_i$ and $\E_{i+2}|\!|W_i$ under $\pi :W^*_i\to
W_i$ by $\E_i|W^*_i$ and $\E_{i+2}|\!|W^*_i$ respectively. Note
that the same $\E_{i+2}|\!| W^*_i$ could have been directly
defined by specializing, similar to case 4 above.

In summary, for a collection of local systems $\E_i$ on $U_i$, 
we have various pullbacks or specializations associated with it follows: 

(i) On $R_i$ we have local systems $\E_i |R_i$ and $\E_{i+1}
|\!| R_i$, for $i\le d-1$

(ii) On $W_i$ we have local systems $\E_i |W_i$ and
$\E_{i+2}|\!|W_i$, for $i\le d-2$

(iii) On $W^*_i$ we have local systems $\E_i|W^*_i$, $\E_{i+1}|\!|W^*_i$, 
and $\E_{i+2}|\!|W^*_i$, for $i\le d-2$ .

\rem On $W^*_i$, we have a canonical identification between 
the sheaf $\E_{i+1}|\!| W^*_i$ (defined as in case (4) above)
and the sheaf $(\E_{i+1}|\!|R_{i+1})|\!|W^*_i$ (defined as in
case (5) above). 

\rem We have earlier defined central elements $\t_i(c)$ in the
fundamental group of each connected component $R_i(c)$ of $R_i$
(see section 2). Given a linear system $\F$ on $R_i$, we denote
by $\t_i$ the automorphism of $\F$ induced by the monodromy
action of the central element $\t_i(c)$ on $R_i(c)$. 

\begin{definition}\rm
A {\bf Verdier object} $(\E_i,C_i,V_i)$ on $(X,Y)$ consists of
the following. 

(1) For each $m\le i\le d$, 
$\E_i$ is a local system on $U_i$ (the ranks of the local systems are
not necessarily constant.) 

(2) For each $m\le i\le d-1$, 
$C_i:(\E_{i+1}|\!| R_i) \to (\E_i| R_i)$ and 
$V_i:(\E_i |R_i) \to (\E_{i+1} |\!| R_i)$ are 
homomorphisms of local systems, such that 

\begin{eqnarray*}
V_iC_i &=& 1-\t_i {\mbox{~{\rm on}~}} 
                     \E_{i+1}|\!| R_i \\ 
C_iV_i &=& 1-\t_i {\mbox{~{\rm on}~}} 
                     \E_i| R_i \\ 
\end{eqnarray*}

(3) Let $m\le i \le d-2$. 
Let $\pi: W^*_i\to W_i$ be the covering projection induced by 
$\pi:Z^*_i\to Z_i$. Let 
$$a_{i+2}:\E_{i+2}|\!|W_i \to \pi_*\pi^*(\E_{i+2}|\!|W_i) 
= \pi_*(\E_{i+2}|\!|W^*_i)$$
$$a_i:\E_i|W_i \to \pi_*\pi^*(\E_i|W_i) = \pi_*(\E_i|W^*_i)$$
be adjunction maps, and let the cokernels of these maps be denoted by 
$$q_{i+2}:\pi_*(\E_{i+2}|\!|W_i)\to Q_{i+2}$$
$$q_i:\pi_*(\E_i|W_i)\to Q_i$$
Then we impose the
requirement that the composite map 
$$\E_{i+2}|\!|W_i\stackrel{a_{i+2}}{\to} \pi_*(\E_{i+2}|\!|W^*_i) 
\stackrel{\pi_*(C_{i+1}|\!|W^*_i)}{\to}  
\pi_*(\E_{i+1}|\!| W^*_i) 
\stackrel{\pi_*(t_i|W^*_i)}{\to} \pi_*(\E_i | W^*_i)
\stackrel{q_i}{\to} Q_i  $$
is zero.

(4) Similarly, we demand that for all $m\le i\le d-2$ the
composite map 
$$Q_{i+2}\stackrel{q_{i+2}}{\leftarrow} \pi_*(\E_{i+2}|\!| W^*_i) 
\stackrel{\pi_*(V_{i+1}|\!|W^*_i)}{\leftarrow} 
\pi_*(\E_{i+1}|\!| W^*_i) 
\stackrel{\pi_*(V_i|W^*_i)}{\leftarrow} \pi_*(\E_i | W^*_i) 
\stackrel{a_i}{\leftarrow} \E_i|W_i$$
is zero.

(5) Note that as $\pi:W^*_i\to W_i$ is a 2-sheeted cover, 
for any sheaf $\F$ on $W_i$ the new sheaf $\pi_*\pi^*(\F)$ on $W_i$ has 
a canonical involution coming from the deck transformation for 
$W^*_i\to W_i$ which transposes the two points over any base point. 
In particular, the local systems
$\pi_*(\E_{i+2}|\!|W^*_i)=\pi_*\pi^*(\E_{i+2}|\!|W_i)$ and 
$\pi_*(\E_i|W^*_i)=\pi_*\pi^*(\E_i|W_i))$ have canonical
involutions, which we denote by $\nu$. We demand that the
following diagram should commute. 

Diagram III.
$$\begin{array}{ccccc}
\pi_*(\E_{i+1} |\!| W^*_i) & \stackrel{\pi_*(V_{i+1}|\!|W^*_i)}{\to} & 
\pi_*(\E_{i+2}|\!|W^*_i) & \stackrel{\nu}{\to} &
\pi_*(\E_{i+2}|\!|W^*_i) \\ 
{\scriptstyle \pi_*(C_{i+1}|\!|W^*_i)}\downarrow & & & & 
\downarrow {\scriptstyle \pi_*(C_i|W^*_i)}\\ 
\pi_*(\E_i|W^*_i) & \stackrel{\nu}{\to} & \pi_*(\E_i|W^*_i)  
& \stackrel{\pi_*(V_i|W^*_i)}{\to} & \pi_*(\E_{i+1} |\!| W^*_i) \\
\end{array}$$

\end{definition}

\rem 
As the adjunction maps are injective (in particular as $a_i$ is
injective), the condition (3) is equivalent to demanding the
existence of a unique $f$ which makes the following diagram commute.

Diagram I.
$$\begin{array}{ccccc}
\E_{i+2}|\!|W_i & & \stackrel{f}{\longrightarrow} & & \E_i|W_i \\
a_{i+2}\downarrow & & & & \downarrow a_i\\
\pi_*(\E_{i+2}|\!| W^*_i) &
\stackrel{\pi_*(C_{i+1}|\!|W^*_i)}{\to} & 
\pi_*(\E_{i+1}|\!| W^*_i) &  
\stackrel{\pi_*(C_i|W^*_i)}{\to} & 
\pi_*(\E_i | W^*_i)\\
\end{array}$$
Similarly, the condition (4) is equivalent to the following:
there must exist a unique homomorphism $g$ 
which makes the following diagram commute.

Diagram II.
$$\begin{array}{ccccc}
\E_{i+2}|\!|W_i & & \stackrel{g}{\longleftarrow} & & \E_i|W_i \\
a_{i+2}\downarrow & & & & \downarrow a_i\\
\pi_*(\E_{i+2}|\!| W^*_i) &
\stackrel{\pi_*(V_{i+1}|\!|W^*_i)}{\leftarrow} & 
\pi_*(\E_{i+1}|\!| W^*_i) &  
\stackrel{\pi_*(V_i|W^*_i)}{\leftarrow} & 
\pi_*(\E_i | W^*_i)\\
\end{array}$$

It can be seen that the above definition of a Verdier object on
$(X,Y)$ reduces in the case of a polydisk to the hypercube
description of perverse sheaf on a polydisk. As
Verdier objects, perverse sheaves, and the specialization
functors are all local in nature, we get the following by gluing up. 

\begin{proposition}
There is an equivalence of categories between the category of \\
Verdier objects and the category of perverse sheaves on $(X,Y)$.
\end{proposition}

\subsection{Moduli for perverse sheaves}

As the various fundamental groups are finitely generated, the 
definition of a Verdier object has an immediate translation
in terms of {\bf quivers}, that is, diagrams of finite
dimensional vector 
spaces and linear maps, by means of remark \ref{composite}.
We define a {\bf family of Verdier objects}
parametrized by some space $T$ as a family of such quivers over
$T$, in which vector spaces are replaced by vector bundles over
$T$ and linear maps (or group representations) are replaced by 
endomorphisms of the bundles. These obviously form an algebraic
stack in the sense of Artin, if we work in the category of
schemes  over $\C$.

When we fix the ranks $n_{i,a}$ of the restrictions of local
systems $\E_i$ on connected components $X_{i,a}$ of $X_i$, and
go modulo the conjugate 
actions of the various $GL(n_a)$ (this is exactly as in the
section 6 of [N-S] so we omit the details), 
we get an affine scheme of finite type over $\C$ 
as the moduli of Verdier objects with given ranks. 
The points of this moduli space are Jordan-Holder classes (that is, 
semisimplifications) of Verdier objects.

The above definition of families and construction of moduli is
independent (upto isomorphism) of the choices of base points and
generators for the various fundamental groups.

We define an {\bf algebraic} (or {\bf holomorphic}) {\bf family
of perverse sheaves on $(X,Y)$} to be an algebraic (or
holomorphic) family  
of Verdier objects, parametrized by a complex scheme (or a
complex analytic space) $T$. Therefore, we have

\begin{proposition}\label{points}
There exists a coarse moduli scheme $\cal P$ for perverse sheaves on
$(X,Y)$ of fixed numerical type. The scheme $\P$ is an affine
scheme of finite type over $\C$, and points of $\P$ correspond
to Jordan-Holder classes of perverse sheaves.
\end{proposition}

\section{The Riemann-Hilbert morphism}
In section 7.1 we define  an analytic
morphism $\rh$ from the 
stack (or moduli) of pre-$\D$-modules to the stack (or moduli)
of Verdier objects, which reperesents the de Rham functor. 
Note that even if both sides are algebraic, the map is only
analytic in general, as it involves integration in order to
associate to a connection its monodromy.

Next (in section 7.2) we prove some properties of the above
Riemann-Hilbert morphism $\rh$, in particular that it is a local
isomorphism at points representing pre-$\D$-modules which have
good residual eigenvalues. This generalizes the rigidity results
in [N] and [N-S]. 

\subsection{Definition of the Riemann-Hilbert morphism}
The following allows us to go from pre-$\D$-modules to 
perverse sheaves. 

\begin{proposition}\label{map} Let $X$ be a nonsingular variety,
$Y\subset X$ a divisor with normal crossing, and let $Y^*\to Y$ the
normalization of $Y$, with $f:Y^*\to X$ the composite map. Let
$N$ be the normal bundle to $f: Y^* \to X$, and let $F\subset
N$ be the closed subset of the total space of $N$ defined by
vestors tangent to branches of $Y$ (in particular, this includes
the zero section $Y^*$ of $N$). Let $\pi : N_f\to Y$ be the bundle
projection.  Then we have 

(1) If $E$ is a vector bundle on $Y^*$ together with the
structure of a $\D_N[\log F]|Y^*$-module, then $\pi^*F$ is naturally
a $\D_N[\log F]$-module.

(2) Let $E$ be a vector bundle on $X$
together with the structure of a $\D_X[\log Y]$-module, and let
$E|Y^*$ be its pullback to $Y^*$, which is naturally a module
over $f^*\D_X[\log Y]$. Then $\pi^*(E|Y^*)$ is naturally a
$\D_N[\log F]$-module. 

(3) If $E$ is as in (2) above and if the residual eigenvalues of
$E$ do not differ by non-zero 
integers on any component of $Y^*$, then the local system
$(\pi^*(E|Y^*))^{\nabla}$ 
on $N-F$ of integrable sections of $\pi^*(E|Y^*)$ is canonically
isomorphic to the specialization of the local system
$(E|X-Y)^{\nabla}$ on $X-Y$ of
integrable sections of $E$.

(4) If $E_T$ is a holomorphic family of vector bundles with
integrable logarithmic connections on $(X,Y)$ parametrized by a
complex analytic space $T$, such that each $E_t$ has good residual
eigenvalues, then the corresponding local systems on $N-F$ given
by (3) form a holomorphic family of local systems on $N_F$
parametrized by $T$.
\end{proposition}

\proof The statement (1) is a special case of the following more
general statement. Let
$S$ be any nonsingular variety, $\pi:N\to S$ any geometric
vector bundle on $S$, 
and $F\subset N$ a normal crossing divisor in the total space
of $N$ such that analytic (or \'etale) locally $F$ is the
union of $r$ vector subbundles of $N$ of rank $r-1$ where $r$ is the
rank of $N_S$. Then for any vector bundle $E$ on $S$ together
with the structure of a $\D_N[\log F]$-module, the vector bundle
$\pi^*(E)$ on $N$ has a natural structure of a $\D_N[\log
F]$-module. This can be seen by choosing analytic local
coordinates $(x_1,\ldots,x_m, y_1, \ldots, y_r)$ on $N$ where 
$(x_i)$ are local coordinates on $S$ and $y_i$ are 
linear coordinates on the fibers such that $F$ is locally
defined by $\prod_iy_i=0$, and defining a logarithmic connection
on $\pi^*(E)$ in terms of the actions of $\partial/\partial x_i$
and $y_i\partial/\partial y_i$ given by
\begin{eqnarray*}
\nabla_{\partial/\partial x_i}(g(y)\otimes_{\O_S} e) 
& = & g(y)\otimes_{\O_S} \nabla_{\partial/\partial x_i} e \\
\nabla_{y_i\partial/\partial y_i}(g(y)\otimes_{\O_S} e) 
& = & (y_i(\partial/\partial y_i)g(y))\otimes_{\O_S}  e 
\end{eqnarray*}

The statement (2) follows from the canonical isomorphism between
$f^*(\D_X[\log Y])$ and $\D_N[\log f]|Y^*$. The statements (3)
and (4) follows over polydisks from the relation between
V-filtrations and specializations, and we have the 
global statements by gluing. 

\rem\label{notalg} As we have to integrate in order
to associate its monodromy representation to an integrable
connection, even if $T$ was associated to an algebraic variety
and the original family $(E_T, \, \nabla_T)$ was algebraic, the
associated family of local systems will in general only be an
analytic family.

\rem It is sometimes erroneously believed that if $E$ is a
a locally free logarithmic connection on $(X,Y)$, and
$Y$ is nonsingular, then the restriction $E|Y$ has a natural
connection on it. This is false in general, and is correct in
some special case if the exact sequencence $0\to \O_Y\to
T_X[\log Y]|Y \to T_Y$ has a natural splitting in some special
case under consideration. 

The above proposition allows us to directly associate a Verdier object
$(\E_i, C_i, V_i)$ to a pre-$\D$-module $(E_i, s_i, t_i)$ which
has good residual eigenvalues in the sense of definition
\ref{goodresieigen}, as follows.

\begin{definition}\label{defrh}\rm
We put $\E_d$ to be the local system on $X-Y$ given by $E_d$,
and for $i\le d-1$ we define
$\E_i$ to be the local system on $U_i$ associated to
the logarithmic connection on $(N_i,F_i)$ associated to $E_i$ by
proposition \ref{map}(1). For $m\le j\le d-1$ we put $C_j$ to be
the map induced (using the statements (2) and (3) in proposition
\ref{map} above) by $\pi_j^*(t_j)$ where
$\pi_j:N_{j,j+1}\to Y^*_j$ is the bundle projection, and define
$V_j$ by the formula 
$$V_j = {\exp(2\pi i\theta_j)-1\over \theta_j} \pi^*(s_j)$$
\end{definition}

It is immediate from the definitions that this gives a Verdier
object starting from a pre-$\D$-module.
By proposition \ref{map}(4), the above association is well
behaved for analytic families, and gives rise to 
an analytic family of Verdier objects when we apply it to an
analytic family of pre-$\D$-modules.  
This gives the {\bf Riemann-Hilbert morphism $\rh$} at the level of
analytic stacks. 

\rem As a consequence of remark \ref{notalg}, 
the above morphism $\rh$ of stacks is holomorphic but not
algebraic.  

\rem By its definition, the Verdier object $\E=(\E_i,C_i,V_i)$
associated by definition \ref{defrh} to a pre-$\D$-module $\EE
=(E_i,t_i,s_i)$ with good residual eigenvalues defines the
perverse sheaf associated by the de Rham functor to the
$\D$-module $M$ associated to $\EE$ in section 4.2 above, as
follows locally from the hypercube descriptions of these objects
restricted to polydisks.

If a semistable pre-$\D$-module has good residual eigenvalues,
then the graded object associated to any S-filtration again has
good residual eigenvalues by remark \ref{goodresi2}. It follows
that the condition that the residual eigenvalues be good defines
an analytic open subset $\M_o$ of the moduli space $\M$ by theorem
\ref{maintheorem}(2). 
It can be proved as in Lemma 7.3 of [N-S] and its following
discussion, using analytic properties of a GIT quotient proved
by Simpson in [S], that our association of a Verdier object to a
pre-$\D$-module with good residual eigenvalues now descends to
an {\bf analytic} morphism $\rh :\M_o\to P$ from $\M_o$ to the
moduli $\P$ of Verdier objects. This is the Riemann-Hilbert
morphism at the level of the moduli spaces.

\subsection{Properties of the Riemann-Hilbert morphism}

This section contains material which is a straightforward
generalization of [N-S], so we omit the details.

In the reverse direction, can construct a pre-$\D$-module with
good residual eigenvalues over a given Verdier object by 
using repeatedly the {\bf Deligne construction}. This gives the
{\bf surjectivity} of 
the Riemann-Hilbert morphism. The exponential map 
$M(n,\C)\to GL(n,\C)$ is a submersion at points where the 
eigenvalues do not differ by $2\pi i$ times a nonzero integer.
Using this, we can extend the Deligne construction to families of
local systems parametrized by Artin local rings (e.g., 
$\C[\epsilon ]/(\epsilon ^2)$) to get families of logarithmic
connections with good residual eigenvalues. From this it follows
that the Riemann-Hilbert morphism is {\bf surjective at tangent
level} at points above having good residual eigenvalues.  

Proposition 5.3 of [N] shows that for a meromorphic connection
$M$ on $X$ with regular singularities on a normal crossing divisor $Y$, 
any locally free logarithmic lattice whose residual eigenvalues
do not differ by positive integers is infinitesimally rigid. In
Proposition 8.6 of [N-S], this is extended to pre-$\D$-modules
on $(X,Y)$ when $Y$ is smooth, by analyzing the derivative of a
map of the form  
$$(s,\,t)\mapsto (s,\,t{e^{st}-1\over st})$$
for matrices $s$ and $t$ (lemma 3.10 of [N-S]). 
By a similar proof applied to the formula given in definition
\ref{defrh}, we have the following when $Y$ is normal crossing.

\begin{proposition} {\rm ({\bf Infinitesimal rigidity}):}   
A pre-$\D$-module on $(X,Y)$ with good residual eigenvalues does
not admit any nontrivial infinitesimal deformations such that the
associated $\D$ module (or perverse sheaf) on $(X,Y)$ is constant.  
\end{proposition}

This shows that the Riemann-Hilbert morphism is a 
{\bf tangent level isomorphism} of stacks at points above with
good residual eigenvalues. 

The above properties, which are for the Riemann-Hilbet morphism
as a morphism of analytic stacks, are valid by \ref{maintheorem}
and \ref{points} for the morphism $\rh : \M_o \to \P$ at the level
of the two moduli spaces at stable points of $\M_o$
which go to simple Verdier objects.

\section*{References} \addcontentsline{toc}{section}{References}
[G-G-M] Galligo, Granger, Maisonobe : 
$\D$-modules et faisceaux pervers dont le support
singulier est un croissement normal. Ann. Inst. Fourier,
Grenoble 35 (1985) 1-48.

[G-M-V] Gelfand, MacPherson, Vilonen : Perverse sheaves and
quivers. Duke Math. J. 83 (1996) 621-643.  

[L] Laumon, G. : Champs alg\'ebriques. Preprint no. 88-33,
Universit\'e Paris Sud, 1988.

[M-V] MacPherson and Vilonen : Elementary construction of
perverse sheaves. Invnt. Math. 84 (1986), 403-435.

[Mal] Malgrange, B. : Extension of holonomic $\D$-modules, in
Algebraic Analysis (dedicated to M. Sato), M. Kashiwara and T.
Kawai eds., Academic Press, 1988.

[Ne] Newstead, P.E. : {\sl Introduction to moduli problems and
orbit spaces}, TIFR lecture notes, Bombay (1978).

[N] Nitsure, N. : Moduli of semistable logarithmic connections.
J. Amer. Math. Soc. 6 (1993) 597-609.

[N-S] Nitsure, N. and Sabbah, C. : Moduli of pre-$\D$-modules,
perverse sheaves, and the Riemann-Hilbert morphism -I, 
Math. Annaln. 306 (1996) 47-73.

[S] Simpson, C. : Moduli of representations of the fundamental
group of a smooth projective variety - I, Publ. Math. I.H.E.S.
79 (1994) 47-129.

[V1] Verdier, J.-L. : Extension of a perverse sheaf across a closed
subspace, \\ Ast\'erisque 130 (1985) 210-217. 

[V2] Verdier, J.-L. : Prolongements de faisceaux pervers
monodromiques, \\
Ast\'erisque 130 (1985) 218-236.

\bigskip

Address:

School of Mathematics, Tata Institute of Fundamental Research,
Homi Bhabha Road, Mumbai 400 005, India. e-mail:
nitsure@math.tifr.res.in

\bigskip

\centerline{09-III-1997}

\end{document}